\documentclass[11pt]{article}

\usepackage{amsfonts}
\usepackage{amscd}
\usepackage{amssymb}
\usepackage{amsthm}
\usepackage{amsmath, xspace}
\usepackage{blkarray}
\usepackage{bm}
\usepackage{cancel}
\usepackage{color}
\usepackage{graphics}
\usepackage{graphicx}
\usepackage{enumitem}
\usepackage{fancyhdr}
\usepackage{mathdots}
\usepackage{mathrsfs}
\usepackage{multicol}
\usepackage{stmaryrd}
\usepackage{tkz-fct}
\usepackage{tkz-euclide}
\usepackage{ytableau}
\usepackage[all]{xy}

\usepackage[plainpages,backref]{hyperref}

\usepackage{lscape}

\theoremstyle{plain}
\newtheorem{thm}{Theorem}[section]

\newtheorem{cor}[thm]{Corollary}
\theoremstyle{definition}

\theoremstyle{example}

\theoremstyle{remark}

\numberwithin{equation}{section}

\def\cL{\mathcal{L}}

\def\CC{\mathbb{C}}

\def\FF{\mathbb{F}}

\def\QQ{\mathbb{Q}}
\def\RR{\mathbb{R}}

\def\ZZ{\mathbb{Z}}

\def\fp{\mathfrak{p}}

\def\Card{\mathrm{Card}}

\usepackage{etex} 
\usepackage{pictexwd}
\usepackage{tikz}
	\usepgflibrary[patterns] 
	\usetikzlibrary{patterns} 
	\usepgflibrary{shapes.geometric}
        \usetikzlibrary{arrows, calc, positioning}

\newcounter{r}
\newcounter{s}

\newcommand\Part[1]{
        \setcounter{r}{1}
	 \foreach \x in {#1}{
 	{\ifnum\value{r}=1
		\draw (0,\value{r}-1)--(\x,\value{r}-1); 
		\fi}
	\draw (0,\value{r}) to (\x,\value{r});
   	\foreach \y in {0, ..., \x} {\draw (\y,\value{r})--(\y,\value{r}-1);}
	\addtocounter{r}{1}
 }}
 
\newcommand\Tableau[1]{
        \foreach \x [count = \c from 1] in {#1} {
		\foreach \y [count = \d from 1] in \x{
			\node at (\d-.5,\c-.5) {\scriptsize$\y$}; 
			\draw (\d,\c) to (\d,\c-1);
			{\ifnum\d=1
				\draw (0,\c) to (0,\c-1);
				\fi}
			\setcounter{r}{\d}
		}
		{\ifnum\c=1
			\draw (0,0)--(\value{r},0);
			\fi}
		\draw(0,\c) to (\value{r},\c);
		\draw(0,\c-1) to (\value{r},\c-1);
		\setcounter{s}{\c}}}
		
\newcommand\sTableau[1]{
        \foreach \x [count = \c from 1] in {#1} {
		\foreach \y [count = \d from 1] in \x{
			\node at (\d-.5,\c-.5) {\tiny$\y$}; 
			\draw (\d,\c) to (\d,\c-1);
			{\ifnum\d=1
				\draw (0,\c) to (0,\c-1);
				\fi}
			\setcounter{r}{\d}
		}
		{\ifnum\c=1
			\draw (0,0)--(\value{r},0);
			\fi}
		\draw(0,\c) to (\value{r},\c);
		\setcounter{s}{\c}}}

\newcommand{\PartB}[1]{
 \foreach \x [count=\s from 1] in {#1}{
 	{\ifnum\s=1
		\draw (0,\s-1)--(\x,\s-1); 
		\fi}
   \draw (0,\s) to (\x,\s);
   \foreach \y in {0, ..., \x} {\draw (\y,\s)--(\y,\s-1);}
 }}

\tikzstyle{V}=[draw, fill =black, circle, inner sep=0pt, minimum size=1.5pt]
\tikzstyle{wV}=[draw, fill =white, circle, inner sep=0pt, minimum size=4.5pt]
\tikzstyle{bV}=[draw, fill =black, circle, inner sep=0pt, minimum size=4.5pt]
\tikzstyle{over}=[draw=white,double=black,line width=2pt, double distance=.5pt]

\def\Over[#1,#2][#3,#4]{ 
	\draw[style=over]   (#2,#1) .. controls ++(#4*.5-#2*.5,0) and ++(-#4*.5+#2*.5,0) .. (#4,#3);}
\def\Under[#1,#2][#3,#4]{ 
	\draw  (#2,#1) .. controls ++(#4*.5-#2*.5,0) and ++(-#4*.5+#2*.5,0) .. (#4,#3);}
\def\Cross[#1,#2][#3,#4]{
	\Under[#3,#2][#1,#4]\Over[#1,#2][#3,#4]}

\def\Tops[#1][#2][#3]{
	\foreach\x in {#1}{
		\draw (#2,\x+.15) -- (#2+.1, \x+.15) (#2, \x-.15) -- (#2+.1, \x-.15) ;
		\draw (#2+.1,\x) arc (0:360:.75mm and 1.5mm);}
	\foreach \x in {1,...,#3} {\draw (#2,\x)  to (#2+.05,\x); \node[V] at (#2+.05,\x){};}
	}
\def\Bottoms[#1][#2][#3]{
	\foreach\x in {#1}{
		\draw (#2, \x+.15) -- (#2-.1, \x+.15) (#2, \x-.15) -- (#2-.1, \x-.15) ;
		\draw (#2-.1, \x+.15) arc (90:270:.75mm and 1.5mm);}
	\foreach \x in {1,...,#3} {\draw (#2, \x)  to (#2-.05, \x); \node[V] at (#2-.05, \x){};}
	}
\def\Caps[#1][#2,#3][#4]{
	\Tops[#1][#3][#4]
	\Bottoms[#1][#2][#4]
	}
\def\Pole[#1][#2,#3]{
	\shade[left color=white,right color=white] (#2,#1+.15) rectangle (#3,#1-.15);
	\draw[over] (#2,#1+.15) to (#3,#1+.15) (#2,#1-.15) to (#3,#1-.15) ;}
\def\Label[#1,#2][#3][#4]{
	\node[right] at (#2+.1,#3) {#4};
	\node[left] at (#1-.1,#3) {#4};		}
\def\Nodes[#1][#2]{
	 \foreach \x in {1,...,#2} {\node[V] at (#1,\x){};	}
	}
\def\PoleCaps[#1][#2,#3]{
	\foreach\x in {#1}{
		\draw (#2,\x+.15) -- (#2-.1,\x+.15) (#2,\x-.15) -- (#2-.1,\x-.15) ;
		\draw (#2-.1,\x+.15) arc (0:-180:1.5mm and .75mm);}
	\foreach\x in {#1}{
		\draw (#3,\x+.15) -- (#3+.1,\x+.15) (#3,\x-.15) -- (#3+.1,\x-.15) ;
		\draw (#3+.1,\x+.15) arc (0:360:1.5mm and .75mm);}
	}
\def\PoleTwist[#1,#2]{
	\foreach \x/\y in {-1/1L, -.7/1R, 0/2L, .3/2R}{\coordinate(T\y) at (#2,\x); \coordinate(B\y) at (#1,\x);}
	\draw[thin] (B1R) .. controls ++(#2*.5-#1*.5-.1,0) and ++(-#2*.5+#1*.5-.1,0) ..  (T2R)
			(B1L)   .. controls ++(#2*.5-#1*.5+.1,0) and ++(-#2*.5+#1*.5+.1,0) ..    (T2L) ;
	\draw[line width=2pt, white]
			(#1,.15)  .. controls +(#2*.5-#1*.5,0) and +(-#2*.5+#1*.5,0) ..   (#2,-.85) ;
	\draw[thin,over] 
		(B2R) .. controls ++(#2*.5-#1*.5+.1,0) and ++(-#2*.5+#1*.5+.1,0) ..  (T1R) 
			(B2L)  .. controls +(#2*.5-#1*.5-.1,0) and +(-#2*.5+#1*.5-.1,0) ..   (T1L) ;
			}

\def\SymPolesCaps[#1,#2][#3]{
	\draw (#1,.3) -- (#1-.1,.3) (#1,.15) -- (#1-.1, .15) ;
	\draw (#1-.1, .3) arc (0:-180:2pt and 1.5pt);
	\draw (#1,#3+.7) -- (#1-.1,#3+.7) (#1,#3+.85) -- (#1-.1,#3+.85) ;
	\draw (#1-.1,#3+.85)  arc (0:-180:2pt and 1.5pt);
	\draw (#2,.3) -- (#2+.1, .3) (#2, .15) -- (#2+.1, .15) ;
	\draw (#2+.1, .3) arc (0:360:2pt and 1.5pt);
	\draw (#2, #3+.7) -- (#2+.1, #3+.7) (#2, #3+.85) -- (#2+.1, #3+.85) ;
	\draw (#2+.1, #3+.85) arc (0:360:2pt and 1.5pt);}

\newcommand{\posleq}[1]{
	\hspace{0.1cm}
	\begin{tikzpicture}
	\draw (-0.8ex, -0.5ex) -- (0.8ex, -0.5ex);
	\draw (-0.8ex, 0.4ex) -- (0.7ex, -0.2ex);
	\draw (-0.8ex, 0.4ex) -- (0.7ex, 1ex);
	\draw (0.4ex,0.4ex) --(1.1ex, 0.4ex);
	\draw (0.75ex,0.75ex) --(0.75ex, 0.05ex);
	\end{tikzpicture}
	\hspace{0.1cm}
	}
\newcommand{\negleq}[1]{
	\hspace{0.1cm}
	\begin{tikzpicture}
	\draw (-0.8ex, -0.5ex) -- (0.8ex, -0.5ex);
	\draw (-0.8ex, 0.4ex) -- (0.7ex, -0.2ex);
	\draw (-0.8ex, 0.4ex) -- (0.7ex, 1ex);
	\draw (0.4ex,0.4ex) --(1.1ex, 0.4ex);
	\end{tikzpicture}
	\hspace{0.1cm}
	}
	
\newcommand{\zeroleq}[1]{
	\hspace{0.1cm}
	\begin{tikzpicture}
	\draw (-0.8ex, -0.5ex) -- (0.8ex, -0.5ex);
	\draw (-0.8ex, 0.4ex) -- (0.7ex, -0.2ex);
	\draw (-0.8ex, 0.4ex) -- (0.7ex, 1ex);
	\draw  (0.75ex,0.4ex) ellipse (0.2ex and 0.35ex);
	\end{tikzpicture}
	\hspace{0.1cm}
	}
	
\newcommand{\posgeq}[1]{
	\hspace{0.1cm}
	\begin{tikzpicture}
	\draw (-0.8ex, -0.5ex) -- (0.8ex, -0.5ex);
	\draw (0.8ex, 0.4ex) -- (-0.7ex, -0.2ex);
	\draw (0.8ex, 0.4ex) -- (-0.7ex, 1ex);
	\draw (-0.4ex,0.4ex) --(-1.1ex, 0.4ex);
	\draw (-0.75ex,0.75ex) --(-0.75ex, 0.05ex);
	\end{tikzpicture}
	\hspace{0.1cm}
	}
\newcommand{\neggeq}[1]{
	\hspace{0.1cm}
	\begin{tikzpicture}
	\draw (-0.8ex, -0.5ex) -- (0.8ex, -0.5ex);
	\draw (0.8ex, 0.4ex) -- (-0.7ex, -0.2ex);
	\draw (0.8ex, 0.4ex) -- (-0.7ex, 1ex);
	\draw (-0.4ex,0.4ex) --(-1.1ex, 0.4ex);
	\end{tikzpicture}
	\hspace{0.1cm}
	}
	
\newcommand{\zerogeq}[1]{
	\hspace{0.1cm}
	\begin{tikzpicture}
	\draw (-0.8ex, -0.5ex) -- (0.8ex, -0.5ex);
	\draw (0.8ex, 0.4ex) -- (-0.7ex, -0.2ex);
	\draw (0.8ex, 0.4ex) -- (-0.7ex, 1ex);
	\draw  (-0.75ex,0.4ex) ellipse (0.2ex and 0.35ex);
	\end{tikzpicture}
	\hspace{0.1cm}
	}

\newcommand{\posl}[1]{
	\hspace{0.1cm}
	\begin{tikzpicture}
	\draw (-0.8ex, 0.4ex) -- (0.7ex, -0.2ex);
	\draw (-0.8ex, 0.4ex) -- (0.7ex, 1ex);
	\draw (0.4ex,0.4ex) --(1.1ex, 0.4ex);
	\draw (0.75ex,0.75ex) --(0.75ex, 0.05ex);
	\end{tikzpicture}
	\hspace{0.1cm}
	}
\newcommand{\negl}[1]{
	\hspace{0.1cm}
	\begin{tikzpicture}
	\draw (-0.8ex, 0.4ex) -- (0.7ex, -0.2ex);
	\draw (-0.8ex, 0.4ex) -- (0.7ex, 1ex);
	\draw (0.4ex,0.4ex) --(1.1ex, 0.4ex);
	\end{tikzpicture}
	\hspace{0.1cm}
	}
	
\newcommand{\zerol}[1]{
	\hspace{0.1cm}
	\begin{tikzpicture}
	\draw (-0.8ex, 0.4ex) -- (0.7ex, -0.2ex);
	\draw (-0.8ex, 0.4ex) -- (0.7ex, 1ex);
	\draw  (0.75ex,0.4ex) ellipse (0.2ex and 0.35ex);
	\end{tikzpicture}
	\hspace{0.1cm}
	}
	
\newcommand{\posg}[1]{
	\hspace{0.1cm}
	\begin{tikzpicture}
	\draw (0.8ex, 0.4ex) -- (-0.7ex, 1ex);
	\draw (0.8ex, 0.4ex) -- (-0.7ex, -0.2ex);
	\draw (-0.4ex,0.4ex) --(-1.1ex, 0.4ex);
	\draw (-0.75ex,0.75ex) --(-0.75ex, 0.05ex);
	\end{tikzpicture}
	\hspace{0.1cm}
	}
\newcommand{\negg}[1]{
	\hspace{0.1cm}
	\begin{tikzpicture}
	\draw (0.8ex, 0.4ex) -- (-0.7ex, -0.2ex);
	\draw (0.8ex, 0.4ex) -- (-0.7ex, 1ex);
	\draw (-0.4ex,0.4ex) --(-1.1ex, 0.4ex);
	\end{tikzpicture}
	\hspace{0.1cm}
	}
	
\newcommand{\zerog}[1]{
	\hspace{0.1cm}
	\begin{tikzpicture}
	\draw (0.8ex, 0.4ex) -- (-0.7ex, -0.2ex);
	\draw (0.8ex, 0.4ex) -- (-0.7ex, 1ex);
	\draw  (-0.75ex,0.4ex) ellipse (0.2ex and 0.35ex);
	\end{tikzpicture}
	\hspace{0.1cm}
	}

\makeatletter
\renewcommand{\@makefnmark}{\mbox{\textsuperscript{}}}
\makeatother

\title{Ian G.~Macdonald: Works of Art}
\author{
Arun Ram\quad\ \ email:\ aram@unimelb.edu.au \\
\\
}
\date{}

\usetikzlibrary{arrows.meta}

\begin{document}

\maketitle


\begin{abstract}
\noindent
Ian Macdonald's works changed our perspective on so many parts of algebraic combinatorics and formal power series. 
This talk will display some selected works of the art of Ian Macdonald, representative of different periods of his oeuvre, 
and analyze how they resonate, both for the past development of our subject and for its future. 
\end{abstract}


\setcounter{tocdepth}{1}
\tableofcontents

\bigskip
\bigskip\noindent
\textbf{Acknowledgments.}
First and foremost my thanks go to Ian Macdonald for his teaching, companionship and for giving me bunches of handwritten notes and copies of his books over the years.  I am so grateful that circumstances were such that I was able to convey these thanks directly to him in person in June 2023.  I thank David Lumsden, Ziheng Zhou, Alex Shields and Dhruv Gupta for energy and insight as we worked through Ian Macdonald's unpublished manuscript on the $n$-line.  I thank Chris Macdonald for reaching out to provide scans and files of the $n$-line manuscript.  I am very grateful to Laura Colmenarejo, Persi Diaconis and Ole Warnaar for helpful suggestions and revisions on this tribute article.

\newpage

\section{Preamble}

This paper was prepared for the occasion of a lecture in tribute to Ian G. Macdonald, delivered at
FPSAC 2024 in Bochum, Germany on 22 July 2024.  I want to express thanks to the Executive Committee of FPSAC,
the Organizing Committee of FPSAC 2024, and to the whole of our FPSAC 2024 community for making this lecture a possibility
and for considering me for its delivery.  Macdonald is my hero, and to be asked to play such a role in his legacy touches me deeply.

\bigskip\bigskip

$$
\begin{matrix}
\vcenter{\hbox{\includegraphics[scale=0.6]{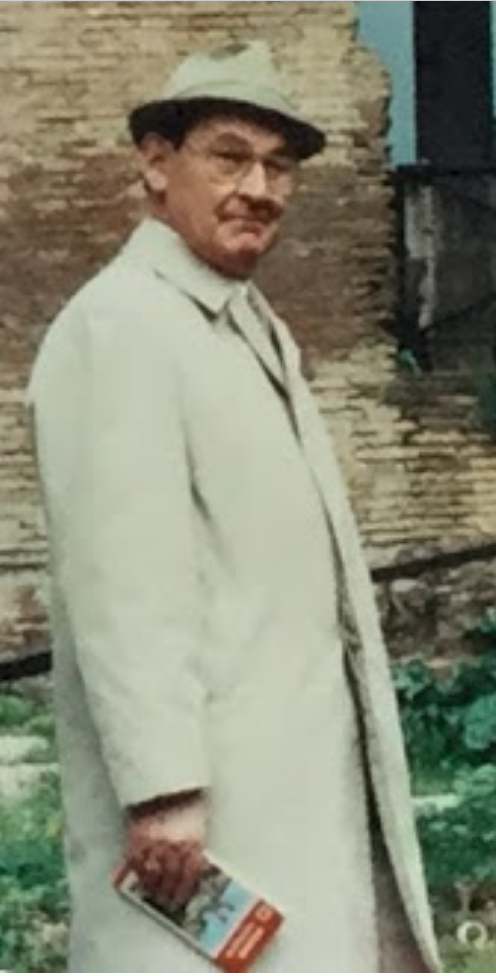}}} \phantom{TTTT}
&\vcenter{\hbox{\includegraphics[scale=0.6]{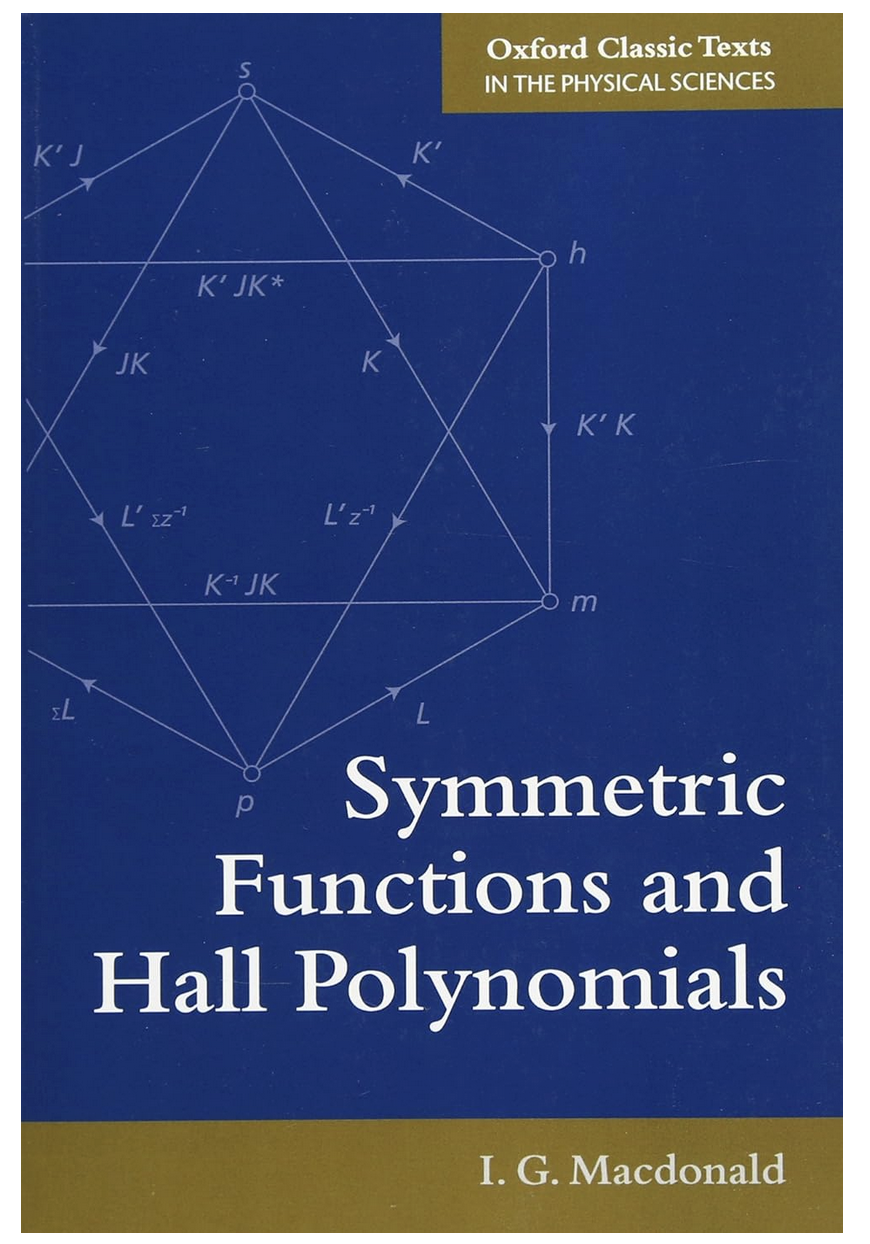}}} \\
\\
\hbox{Ian G. Macdonald} 
&\hbox{The Symmetric Functions Bible}
\end{matrix}
$$
The image of Ian G. Macdonald is from \href{https://sites.google.com/view/garsiafest/mementos}{https://sites.google.com/view/garsiafest/mementos}

\newpage

\section{Tableaux and Macdonald polynomials}

One of our favorite formulas is the formula for the Schur polynomial as a sum over semistandard Young tableaux (SSYTs),
$$s_\lambda = \sum_{T\in B(\lambda)} x^T,
\qquad\hbox{where}\qquad
\begin{array}{c}
B(\lambda) = \{ \hbox{SSYTs of shape $\lambda$}\} \\
\hbox{and} \\
x^T=x_1^{\mathrm{(\#1s\ in\ }T)}\cdots x_n^{\mathrm{(\#}n\mathrm{s\ in\ }T)}.
\end{array}
$$
It is most amazing that if $\delta = (n-1, \ldots, 2,1,0)$ then
$$s_\lambda = \frac{a_{\lambda+\delta}}{a_\delta},
\qquad\hbox{where}\qquad
a_\mu = \sum_{w\in S_n} (-1)^{\ell(w)} w x^\mu
$$
with $x^\mu = x_1^{\mu_1}\cdots x_n^{\mu_n}$ if $\mu  = (\mu_1, \ldots, \mu_n)$.  This second formula for the Schur polynomial is the 
``Weyl character formula'',
which (in this type A case) was one of the first definitions of the Schur function (Jacobi 1841, according to Macdonald).

Macdonald pointed out something spectacular.  The first formula for the Schur polynomial is the special case $q=t$ of the formula
$$P_\lambda(q,t) = \sum_{T\in B(\lambda)} x^T\psi_T(q,t),
\qquad\hbox{where $\psi_T(q,t)$ is given by \eqref{psidef} below,}
$$
and the second formula for the Schur polynomial is a special case of
$$P_\lambda(q,qt) = \frac{A_{\lambda+\delta}(q,t)}{A_\delta(q,t)},
\qquad\hbox{where}\quad
A_{\mu}(q,t) = \sum_{w\in S_n} (-t^{-\frac12})^{\ell(w)} T_w E_{\mu}(q,t)
$$
with $T_w$ and $E_\mu(q,t)$ as defined \eqref{Tdef} and \eqref{Edef} below.  Maybe we think Schur polynomials are cool,
but the Macdonald polynomials $P_\lambda(q,t)$ are two parameters cooler.

\bigskip\bigskip

$$T= \begin{matrix}
\begin{tikzpicture} \Tableau{{3,4},{2,2,4,4},{1,1,1,2,3}} \end{tikzpicture}
\end{matrix}
\qquad\qquad
\begin{matrix}
\begin{tikzpicture}[yscale=-1]
\draw (0,0) -- (5,0) -- (5,1) -- (4,1) -- (4,3)  -- (2,3) -- (2,4) -- (1,4) -- (1,5) -- (0,5) -- (0,0);
\draw[<->] (2.5,0) -- node[anchor=west]{$\scriptstyle{\mathrm{coleg}_\lambda(b)}$} (2.5,1.3) ;
\draw[<->] (2.5,1.7) -- node[anchor=west]{$\scriptstyle{\mathrm{leg}_\lambda(b)}$} (2.5,3) ;
\draw[<->] (0,1.5) -- node[anchor=south]{$\scriptstyle{\mathrm{coarm}_\lambda(b)}$} (2.3,1.5) ;
\draw[<->] (2.7,1.5) -- node[anchor=south]{$\scriptstyle{\mathrm{arm}_\lambda(b)}$} (4,1.5) ;
\draw (2.5,1.5)  node[shape=rectangle,draw]{$b$};
\end{tikzpicture}
\end{matrix}
$$
$$P_\lambda(q,t) = \sum_{T\in B(\lambda)} x^T\psi_T(q,t).$$

%

\newpage

\subsection{$(q,t)$-hooks and the bosonic Macdonald polynomials $P_\lambda(q,t)$}


Let $\lambda\in \ZZ^n_{\ge 0}$ with
$\lambda_1\ge \cdots \ge \lambda_n$ so that $\lambda$ is a partition of length at most $n$.
$$\hbox{A \emph{box in $\lambda$} is a 
pair $b=(r,c)$ with $r\in \{1, \dots, n\}$ and $c\in \{1, \ldots, \lambda_r\}$.}
$$
Identify $\lambda$ with its set of boxes so that
$$\lambda = \{ (r,c)\ |\ \hbox{$(r,c)$ is a box in $\lambda$}\}.$$
For a box $b=(r,c)$ in $\lambda$ define
\begin{align*}
\mathrm{arm}_\lambda(b) &= \mathrm{arm}_\lambda(r,c) =\{ (r,c')\in \lambda\ |\ c'>c\}
\qquad\hbox{and} \\
\mathrm{leg}_\lambda(b) &= \mathrm{leg}_\lambda(r,c) =\{ (r',c)\in \lambda\ |\ r'>r\}.
\end{align*}
A \emph{SSYT (semistandard Young tableau) of shape $\lambda$ filled from $\{1, \ldots, n\}$} is a function
$$T\colon \lambda \to \{1, \ldots, n\}\qquad\hbox{such that} $$
\begin{enumerate}
\item[(a)]If $(r,c), (r+1,c)\in \lambda$ then $T(r,c)< T(r,c+1)$, \quad 
$\begin{matrix}
\begin{tikzpicture}[xscale=0.5,yscale=-0.5]
\draw (0,0) -- (2,0) -- (2,1) -- (0,1) -- (0,0);
\draw (1,0) -- (1,1);
\node at (1,0.5) {$<$};
\end{tikzpicture}
\end{matrix}
$
\item[(b)] If $(r,c), (r,c+1)\in \lambda$ then $T(r,c)\le T(r,c+1)$.
\qquad
$\begin{matrix}
\begin{tikzpicture}[xscale=0.5,yscale=-0.5]
\draw (0,0) -- (1,0) -- (1,2) -- (0,2) -- (0,0);
\draw (0,1) -- (1,1);
\node at (0.5,1) {\rotatebox{90}{$\ge$}};
\end{tikzpicture}
\end{matrix}
$
\end{enumerate}

Let
$$B(\lambda) = \{
\hbox{SSYTs of shape $\lambda$ filled from $\{1, \ldots, n\}$} \}.
$$
Let $T\in B(\lambda)$ and let $b\in \lambda$.  
Let $T(b)$ denote the entry in box $b$ of $T$.
Let $i\in \{1, \ldots, n\}$ with $i> T(b)$.
Define the \emph{$i$-restricted arm length}, \emph{$i$-restricted leg length}, and \emph{$i$-restricted $(q,t)$-hook length} by
\begin{equation}
\begin{array}{l}
a(b, < i) = \Card\{ b'\in \mathrm{arm}_\lambda(b)\ |\ T(b') < i\}, \\
\\
l(b, < i) = \Card\{ b'\in \mathrm{leg}_\lambda(b)\ |\ T(b') < i\},
\end{array}
\qquad\hbox{and}\qquad
h_T(b,< i) = \frac{1-t\cdot q^{a(b,< i)}t^{l(b,< i)} }{1- q\cdot q^{a(b,< i)}t^{l(b, < i)} }.
\label{leihooks}
\end{equation}
For a column strict tableau $T\in B(\lambda)$ define
\begin{equation}
\psi_T(q,t) = \prod_{b\in \lambda} \psi_T(b),
\qquad\hbox{where}\qquad
\psi_T(b) = \prod_{i>T(b), i\in T(\mathrm{arm}_\lambda(b))\atop i\not\in 
T(\mathrm{leg}\lambda(b))}
\frac{ h_T(b, < i)}{h_T(b, < i+1)}.
\label{psidef}
\end{equation}

The \emph{bosonic Macdonald polynomial} is $P_\lambda(q,t)\in \CC[x_1,\ldots, x_n]$ given by
\begin{equation*}
P_\lambda(q,t) = \sum_{T\in B(\lambda)}   x^T \psi_T(q,t), \qquad\hbox{where}\quad x^T = 
x_1^{\mathrm{(\#1s\ in\ }T)}\cdots x_n^{\mathrm{(\#}n\mathrm{s\ in\ }T)}.
\end{equation*}
The 
\emph{Schur polynomial} is
$$s_\lambda = P_\lambda(t,t) = P_\lambda(0,0) = \sum_{T\in B(\lambda)} x^T.$$

\subsection{Electronic and fermionic Macdonald polynomials}

For $i\in \{1, \ldots, n-1\}$ and $f\in \CC[x_1, \ldots, x_n]$ define
$$(s_if)(x_1, \ldots, x_n) = f(x_1, \ldots, x_{i-1}, x_{i+1}, x_i, x_{i+2}, \ldots, x_n)
$$
and
$$
T_if = -t^{-\frac12}f + (1+s_i)\frac{t^{-\frac12}-t^{\frac12}x_i^{-1}x_{i+1}}{1-x_i^{-1}x_{i+1}} f.
$$
If $w\in S_n$ and $w=s_{i_1}\cdots s_{i_\ell}$ is a reduced word for $w$ as a product of $s_i$s then write
\begin{equation}
T_w = T_{i_1}\cdots T_{i_\ell}
\qquad\hbox{and}\qquad
w = s_{i_1}\cdots s_{i_\ell}.
\label{Tdef}
\end{equation}
For $i\in \{1 \ldots, n-1\}$ let
$$\partial_i = (1+s_i)\frac{1}{x_i-x_{i+1}} .$$
The \emph{electronic Macdonald polynomial $E_\mu=E_\mu(q,t)$} is recursively determined by
\begin{enumerate}
\item[(E0)] $E_{(0,\ldots, 0)} = 1$,
\item[(E1)] $E_{(\mu_n+1, \mu_1, \ldots, \mu_{n-1})} = q^{\mu_n}x_1 E_\mu(x_2, \ldots, x_n, q^{-1}x_1),$
\item[(E2)] If $(\mu_1, \ldots, \mu_n)\in \ZZ^n_{\ge 0}$ and $\mu_i>\mu_{i+1}$ then
\begin{equation}
E_{s_i\mu} = \Big( \partial_ix_i - t x_i\partial_i + \frac{(1-t)q^{\mu_i-\mu_{i+1}} t^{v_\mu(i)-v_\mu(i+1)}}
{1-q^{\mu_i-\mu_{i+1}}t^{v_\mu(i)-v_\mu(i+1)} }\Big)E_\mu,
\label{Edef}
\end{equation}
where $v_\mu\in S_n$ is the minimal length permutation such that $v_\mu\mu$ is weakly increasing.
\end{enumerate}
The monomial $x^\mu$ is $x^\mu = x_1^{\mu_1}\cdots x_n^{\mu_n}$.
The world of Macdonald polynomials replaces the monomials $x^\mu$ with electronic Macdonald polynomials $E_\mu$ and
replaces the action of permutations $w$ by the operators $T_w$.

Let $\delta = (n-1, n-2, \ldots, 2,1,0)$ and let $\lambda = (\lambda_1, \ldots, \lambda_n)\in \ZZ^n$ with $\lambda_1\ge \cdots \ge \lambda_n$.
\hfil\break
Let $w_0$ be the longest element of $S_n$ so that $\ell(w_0) = \binom{n}{2}$. Then define
\begin{equation*}
A_{\lambda+\delta}(q,t) = \sum_{w\in S_n} (-t^{-\frac12})^{\ell(w)-\ell(w_0)} T_w E_{\lambda+\delta}(q,t)
\qquad\hbox{and}\qquad
a_{\lambda+\delta} = \sum_{w\in S_n} (-1)^{\ell(w)-\ell(w_0)} w x^{\lambda+\delta}.
\end{equation*}
The $A_{\lambda+\delta}(q,t)$ are the \emph{fermionic Macdonald polynomials} (see \cite[Intro]{CR22} for motivation for the 
`electronic', `bosonic', `fermionic' terminology).

\subsection{The Weyl character formula}

The ``Weyl character formula'' in the next theorem
gives a formula for the bosonic Macdonald polynomial as a quotient of two fermionic Macdonald polynomials.  When $q=t=0$ this
formula becomes the formula for the Schur function as a quotient of two determinants.

\begin{thm} Let $\lambda = (\lambda_1, \ldots, \lambda_n)\in \ZZ_{\ge 0}^n$ with $\lambda_1\ge \cdots \ge \lambda_n$.
\item[(a)]
$$P_\lambda(q,qt) = \frac{A_{\lambda+\delta}(q,t)}{A_\delta(q,t)}
\qquad\hbox{and}\qquad
s_\lambda = \frac{a_{\lambda+\delta}}{a_\delta}.
$$
\item[(b)] 
$$A_\delta(q,t) = \prod_{i<j} (x_i-tx_j)
\qquad\hbox{and}\qquad
a_\delta = \prod_{i<j} (x_i-x_j).
$$
\end{thm}

\newpage

\section{Can you do type B?}

Having worked something out for type A, a natural next problem for our community is to work it out for type B.
Here Macdonald has something interesting to say.
$$\hbox{Which type B?}$$
Because, as Macdonald worked out in his 1972 paper on affine root systems,
$$\hbox{there are 9 different type Bs.}$$
A diagram showing these is given in Section \ref{typeBposet}.

But, there is something wonderful here.  The type $(C^\vee,C)$ root system is one of the type Bs and
$$\hbox{all other type Bs are obtained by specializations from type $(C^\vee,C)$.}$$
This means that, if one wants to compute Macdonald polynomials for any one of the 9 different type Bs, then all one
has to do, is compute the Macdonald polynomials for type $(C^\vee,C)$ and then specialize parameters as appropriate.

Each of the affine root systems of classical type is a subset of the $\ZZ$-vector space spanned by symbols $\varepsilon_1, \ldots, \varepsilon_n$ and $\frac12 \delta$,
$$V_\ZZ = \hbox{$\ZZ$-span}\{ \varepsilon_1, \ldots, \varepsilon_n, \hbox{$\frac12$} \delta\}.
$$
The affine Weyl group $W$ is the group of $\ZZ$-linear transformations of $V_\ZZ$
generated by the transformations $s_0, s_1, \ldots, s_n$ given by:
for $\lambda = \lambda_1 \varepsilon_1 + \cdots + \lambda_n \varepsilon_n + \frac{k}{2}\delta$,
\begin{align}
s_0 \lambda &= - \lambda_1\varepsilon_1 + \lambda_2\varepsilon_2
+\cdots+\lambda_n \varepsilon_n + \big(\hbox{$\frac{k}{2}$}+\lambda_1) \delta,
\nonumber  \\
s_n\lambda &= \lambda_1\varepsilon_1+\cdots +\lambda_{n-1}\varepsilon_{n-1}
-\lambda_n \varepsilon_n + \hbox{$\frac{k}{2}$}\delta, \qquad\hbox{and}
\nonumber
\\
s_i\lambda &= \lambda_1\varepsilon_1 + \cdots + \lambda_{i-1}\varepsilon_{i-1}
+\lambda_{i+1}\varepsilon_i+\lambda_i \varepsilon_{i+1} + \lambda_{i+2}\varepsilon_{i+2}+\cdots +
\lambda_n\varepsilon_n + \hbox{$\frac{k}{2}$} \delta,
\nonumber
\end{align}
for $i\in \{1, \ldots, n-1\}$.  
Each of the affine root systems of classical type is
defined by which orbits of the affine Weyl group $W$ that it contains.  Let
\begin{align*}
O_1 
&= W\cdot\alpha_n = W\cdot \varepsilon_n
= \{ \pm \varepsilon_i+r\delta \ |\  i\in \{1, \ldots, n\}, r\in \ZZ\}, \\
O_2 &= W\cdot 2\alpha_n = W\cdot 2\varepsilon_n
= \{ \pm 2\varepsilon_i + 2r\delta \ |\  i\in \{1, \ldots, n\}, r\in \ZZ\} , \\
O_3 
&= W\cdot \alpha_0 = W\cdot  (-\varepsilon_1+\hbox{$\frac12$}\delta)  
= \{ \pm(\varepsilon_i + \hbox{$\frac12$}(2r+1)\delta \ |\  i\in \{1, \ldots, n\}, r\in \ZZ\}, \\
O_4 
&= W\cdot 2\alpha_0 = W\cdot (-2\varepsilon_1+\delta)
= \{ \pm 2\varepsilon_i + (2r+1)\delta \ |\  i\in \{1, \ldots, n\}, r\in \ZZ\} , \\
O_5  
&= W\cdot \alpha_1 = W\cdot (\varepsilon_1-\varepsilon_2)
= 
\left\{ \begin{array}{l}
\pm(\varepsilon_i + \varepsilon_j)+r\delta  \\
\pm(\varepsilon_i - \varepsilon_j)+r\delta  
\end{array}\ \Big\vert\ 
i,j\in \{1, \ldots, n\}, i<j, r\in \ZZ \right\},
\end{align*}
where
$$
\begin{matrix}
\begin{tikzpicture}[every node/.style={inner sep=1}, scale=1]
	\node at (.5,0) {\tikz \draw[black, ->, >=implies, double distance=2] (0,0) -- +(-.127,0);};
	\node at (5.5,0) {\tikz \draw[black, ->, >=implies, double distance=2] (0,0) -- +(.127,0);};
	\foreach \x in {0,5}{\foreach \y in {-1,1}{\draw (\x,\y*.05) to (\x+1,\y*.05);}}
	\node[wV, label=below left:{$2\alpha_0 = -2\varepsilon_1+\delta$}, label=above left:{$\alpha_0 = -\varepsilon_1+\frac12 \delta$}] (0) at (0,0) {};
	\node[wV, label=below right:{$2\alpha_n = 2\varepsilon_n$}, label=above right:{$\alpha_n = \varepsilon_n$}] (6) at (6,0) {};
	\foreach \x/\y in {1/1,2/2,4/n\!-\!2,5/n\!-\!1}{
		\node[wV] (\x) at (\x,0) {};}
	\node[wV,label=above right:{$\alpha_i = \varepsilon_i - \varepsilon_{i+1}$}] (2) at (2,0) {};
	\draw (1) to (2) (4) to (5);
	\draw[dashed] (2) to (4);
	\end{tikzpicture}\end{matrix}.
$$
With these notations the irreducible affine root systems of classical type (and the appropriate specializations for obtaining the 
Macdonald polynomials of each type from the Macdonald polynomials of type $(C^\vee,C)$) are given by the following diagram.
The middle notation for each root system  is the notation in Macdonald \cite[\S 1.3]{Mac03}, the right notation
is that of Bruhat-Tits \cite{BT72} and the leftt notation is that of Kac \cite[Ch.\ 6]{Kac}.

\begin{landscape}
\subsection{The poset of  affine root systems of classical type} \label{typeBposet}

$$\begin{tikzpicture}[xscale=1.7,yscale=-1.4]
\node (CvC) at (0,0) 
{$\begin{matrix}
(C^\vee_n,C_n)=\hbox{C-BC}_n^{\mathrm{II}}
\\
\begin{matrix}
\begin{tikzpicture}[every node/.style={inner sep=1}, scale=.75]
	\node at (.5,0) {\tikz \draw[black, ->, >=implies, double distance=2] (0,0) -- +(.127,0);};
	\node at (5.5,0) {\tikz \draw[black, ->, >=implies, double distance=2] (0,0) -- +(.127,0);};
	\foreach \x in {0,5}{\foreach \y in {-1,1}{\draw (\x,\y*.05) to (\x+1,\y*.05);}}
	\node[wV, label=below:{\tiny$O_4$}, label=above:{\tiny$O_3$}] (0) at (0,0) {};
	\node[wV, label=below:{\tiny$O_2$}, label=above:{\tiny$O_1$}] (6) at (6,0) {};
	\foreach \x/\y in {1/1,2/2,4/n\!-\!2,5/n\!-\!1}{
		\node[wV,label=above:{\tiny$O_5$}] (\x) at (\x,0) {};}
	\draw (1) to (2) (4) to (5);
	\draw[dashed] (2) to (4);
	\end{tikzpicture}\end{matrix}
\\
P_\mu(x;q,t_n^{1/2},u_n^{1/2},t_0^{1/2},u_0^{1/2},t)
\end{matrix}
$};
\node (CvBC) at (-2,2) {
$\begin{matrix}
(C^\vee_n,BC_n)=\hbox{C-BC}_n^{\mathrm{I}}
\\
\begin{matrix}
\begin{tikzpicture}[every node/.style={inner sep=1}, scale=.75]
	\node at (.5,0) {\tikz \draw[black, ->, >=implies, double distance=2] (0,0) -- +(-.127,0);};
	\node at (5.5,0) {\tikz \draw[black, ->, >=implies, double distance=2] (0,0) -- +(.127,0);};
	\foreach \x in {0,5}{\foreach \y in {-1,1}{\draw (\x,\y*.05) to (\x+1,\y*.05);}}
	\node[wV, label=above:{\tiny$O_3$}] (0) at (0,0) {};
	\node[wV, label=below:{\tiny$O_2$}, label=above:{\tiny$O_1$}] (6) at (6,0) {};
	\foreach \x/\y in {1/1,2/2,4/n\!-\!2,5/n\!-\!1}{
		\node[wV,label=above:{\tiny$O_5$}] (\x) at (\x,0) {};}
	\draw (1) to (2) (4) to (5);
	\draw[dashed] (2) to (4);
	\end{tikzpicture}\end{matrix}
\\
P_\mu(x;q,t_n^{1/2},u_n^{1/2},t_0^{1/2},t_0^{1/2},t)
\end{matrix}
$};
\node (BCC) at (2,2) {
$\begin{matrix}
(BC_n,C_n)=\hbox{C-BC}_n^{\mathrm{IV}}
\\
\begin{matrix}
\begin{tikzpicture}[every node/.style={inner sep=1}, scale=.75]
	\node at (.5,0) {\tikz \draw[black, ->, >=implies, double distance=2] (0,0) -- +(.127,0);};
	\node at (5.5,0) {\tikz \draw[black, ->, >=implies, double distance=2] (0,0) -- +(.127,0);};
	\foreach \x in {0,5}{\foreach \y in {-1,1}{\draw (\x,\y*.05) to (\x+1,\y*.05);}}
	\node[wV, label=below:{\tiny$O_4$}] (0) at (0,0) {};
	\node[wV, label=below:{\tiny$O_2$}, label=above:{\tiny$O_1$}] (6) at (6,0) {};
	\foreach \x/\y in {1/1,2/2,4/n\!-\!2,5/n\!-\!1}{
		\node[wV,label=above:{\tiny$O_5$}] (\x) at (\x,0) {};}
	\draw (1) to (2) (4) to (5);
	\draw[dashed] (2) to (4);
	\end{tikzpicture}\end{matrix}
\\
P_\mu(x;q,t_n^{1/2},u_n^{1/2},t_0^{1/2},1,t)
\end{matrix}
$};
\node (Cv) at (-5,4) {
$\begin{matrix}
D_{n+1}^{(2)} = C_n^\vee = \hbox{C-B}_n
\\
\begin{matrix}
\begin{tikzpicture}[every node/.style={inner sep=1}, scale=.75]
	\node at (.5,0) {\tikz \draw[black, ->, >=implies, double distance=2] (0,0) -- +(-.127,0);};
	\node at (5.5,0) {\tikz \draw[black, ->, >=implies, double distance=2] (0,0) -- +(.127,0);};
	\foreach \x in {0,5}{\foreach \y in {-1,1}{\draw (\x,\y*.05) to (\x+1,\y*.05);}}
	\node[wV, label=above:{\tiny$O_3$}] (0) at (0,0) {};
	\node[wV, label=above:{\tiny$O_1$}] (6) at (6,0) {};
	\foreach \x/\y in {1/1,2/2,4/n\!-\!2,5/n\!-\!1}{
		\node[wV,label=above:{\tiny$O_5$}] (\x) at (\x,0) {};}
	\draw (1) to (2) (4) to (5);
	\draw[dashed] (2) to (4);
	\end{tikzpicture}\end{matrix}
\\
P_\mu(x;q,t_n^{\frac12},t_n^{1/2},t_0^{1/2},t_0^{1/2},t)
\end{matrix}
$};
\node (BBv) at (-2,4) {
$\begin{matrix}
(B_n,B_n^\vee)= \hbox{B-BC}_n
\\
	\begin{matrix}\begin{tikzpicture}[every node/.style={inner sep=1}, scale=.75]
	\node at (6.5,0) {\tikz \draw[black, ->, >=implies, double distance=2] (0,0) -- +(.127,0);};
		\foreach \y in {-1,1}{\draw (6,\y*.05) to (6+1,\y*.05);}
		\node[wV, label=below:{\tiny$O_2$}, label=above:{\tiny$O_1$}] (7) at (7,0) {};
	\foreach \x/\y in {2/2,3/3,5/n\!-\!2,6/n\!-\!1}{
		\node[wV, label=above:{\tiny$O_5$}] (\x) at (\x,0) {};}
	\node[wV, label=above:{\tiny$O_5$}] (0) at (1,-.65) {};
	\node[wV, label=below:{\tiny$O_5$}] (1) at (1,.65) {};
	\draw (1) to (2) (0) to (2) (2) to (3) (5) to (6);
	\draw[dashed] (3) to (5);
	\end{tikzpicture}\end{matrix}
\\
P_\mu(x;q,t_n^{1/2},u_n^{1/2},1,1,t)
\end{matrix}
$};
\node (BC) at (2,4) {
$\begin{matrix}
A_{2n}^{(2)}=BC_n=\hbox{C-BC}_n^{\mathrm{III}}
\\
\begin{matrix}
\begin{tikzpicture}[every node/.style={inner sep=1}, scale=.75]
	\node at (.5,0) {\tikz \draw[black, ->, >=implies, double distance=2] (0,0) -- +(.127,0);};
	\node at (5.5,0) {\tikz \draw[black, ->, >=implies, double distance=2] (0,0) -- +(.127,0);};
	\foreach \x in {0,5}{\foreach \y in {-1,1}{\draw (\x,\y*.05) to (\x+1,\y*.05);}}
	\node[wV, label=below:{\tiny$O_4$}] (0) at (0,0) {};
	\node[wV, label=above:{\tiny$O_1$}] (6) at (6,0) {};
	\foreach \x/\y in {1/1,2/2,4/n\!-\!2,5/n\!-\!1}{
		\node[wV,label=above:{\tiny$O_5$}] (\x) at (\x,0) {};}
	\draw (1) to (2) (4) to (5);
	\draw[dashed] (2) to (4);
	\end{tikzpicture}\end{matrix}
\\
P_\mu(x;q,t_n^{1/2},t_n^{1/2},t_0^{1/2},1,t)
\end{matrix}
$};
\node (C) at (5,4) {
$\begin{matrix}
C_n^{(1)}=C_n=C_n
\\
\begin{matrix}
\begin{tikzpicture}[every node/.style={inner sep=1}, scale=.75]
	\node at (.5,0) {\tikz \draw[black, ->, >=implies, double distance=2] (0,0) -- +(.127,0);};
	\node at (5.5,0) {\tikz \draw[black, ->, >=implies, double distance=2] (0,0) -- +(-.127,0);};
	\foreach \x in {0,5}{\foreach \y in {-1,1}{\draw (\x,\y*.05) to (\x+1,\y*.05);}}
	\node[wV, label=below:{\tiny$O_4$}] (0) at (0,0) {};
	\node[wV, label=below:{\tiny$O_2$}] (6) at (6,0) {};
	\foreach \x/\y in {1/1,2/2,4/n\!-\!2,5/n\!-\!1}{
		\node[wV,label=above:{\tiny$O_5$}] (\x) at (\x,0) {};}
	\draw (1) to (2) (4) to (5);
	\draw[dashed] (2) to (4);
	\end{tikzpicture}\end{matrix}
\\
P_\mu(x;q,t_n^{1/2},1,t_0^{\frac12},1,t)
\end{matrix}
$};
\node (B) at (-2,6) {
$\begin{matrix}
B_n^{(1)} = B_n = B_n
\\
	\begin{matrix}\begin{tikzpicture}[every node/.style={inner sep=1}, scale=.75]
	\node at (6.5,0) {\tikz \draw[black, ->, >=implies, double distance=2] (0,0) -- +(.127,0);};
		\foreach \y in {-1,1}{\draw (6,\y*.05) to (6+1,\y*.05);}
		\node[wV, label=above:{\tiny$O_1$}] (7) at (7,0) {};
	\foreach \x/\y in {2/2,3/3,5/n\!-\!2,6/n\!-\!1}{
		\node[wV, label=above:{\tiny$O_5$}] (\x) at (\x,0) {};}
	\node[wV, label=above:{\tiny$O_5$}] (0) at (1,-.65) {};
	\node[wV, label=below:{\tiny$O_5$}] (1) at (1,.65) {};
	\draw (1) to (2) (0) to (2) (2) to (3) (5) to (6);
	\draw[dashed] (3) to (5);
	\end{tikzpicture}\end{matrix}
\\
P_\mu(x;q,t_n^{1/2},t_n^{1/2},1,1,t)
\end{matrix}
$};
\node (Bv) at (2,6) {
$\begin{matrix}
A_{2n-1}^{(2)} = B_n^\vee =  \hbox{B-C}_n
\\
	\begin{matrix}\begin{tikzpicture}[every node/.style={inner sep=1}, scale=.75]
	\node at (6.5,0) {\tikz \draw[black, ->, >=implies, double distance=2] (0,0) -- +(-.127,0);};
		\foreach \y in {-1,1}{\draw (6,\y*.05) to (6+1,\y*.05);}
		\node[wV, label=below:{\tiny$O_2$}] (7) at (7,0) {};
	\foreach \x/\y in {2/2,3/3,5/n\!-\!2,6/n\!-\!1}{
		\node[wV, label=above:{\tiny$O_5$}] (\x) at (\x,0) {};}
	\node[wV, label=above:{\tiny$O_5$}] (0) at (1,-.65) {};
	\node[wV, label=below:{\tiny$O_5$}] (1) at (1,.65) {};
	\draw (1) to (2) (0) to (2) (2) to (3) (5) to (6);
	\draw[dashed] (3) to (5);
	\end{tikzpicture}\end{matrix}
\\
P_\mu(x;q,t_n^{1/2},1,1,1,t)
\end{matrix}
$};
\node (D) at (0,8) {
$\begin{matrix}
D_n^{(1)}=D_n=D_n
\\
	\begin{matrix}\begin{tikzpicture}[every node/.style={inner sep=1}, scale=.75]
	\foreach \x/\y in {2/2,3/3,5/n\!-\!3,6/n\!-\!2\phantom{-} }{
		\node[wV, label=above:{\tiny$O_5$}] (\x) at (\x,0) {};}
	\node[wV, label=above:{\tiny$O_5$}] (0) at (1,-.5) {};
	\node[wV, label=below:{\tiny$O_5$}] (1) at (1,.5) {};
	\node[wV, label=above:{\tiny$O_5$}] (7) at (7,-.5) {};
	\node[wV, label=below:{\tiny$O_5$}] (8) at (7,.5) {};
	\draw (1) to (2) (0) to (2) (2) to (3) (5) to (6) (6) to (7) (6) to (8);
	\draw[dashed] (3) to (5);
	\end{tikzpicture}\end{matrix}
\\
P_\mu(x;q,1,1,1,1,t)
\end{matrix}
$};
\node (A) at (0,9.5) {$GL_n$};
\node (O) at (0,10.3) {$\emptyset$};
\draw (CvC) -- (CvBC);
\draw (CvC) -- (BCC);
\draw (CvBC) -- (Cv);
\draw (CvBC) -- (BBv);
\draw (BCC) -- (BC);
\draw (BCC) -- (C);
\draw (Cv) -- (B);
\draw (BBv) -- (B);
\draw (BBv) -- (Bv);
\draw (BC) -- (B);
\draw (C) -- (Bv);
\draw (B) -- (D);
\draw (Bv) -- (D);
\draw (D) -- (A);
\draw (A) -- (O);
\end{tikzpicture}
$$
\end{landscape}

%
%

\newpage
\section{Circles and Lines}


Though I don't travel often to England, whenever a trip did bring me to England
I liked to try to stop in and visit Ian and his wife Greta if I could manage it.  
Greta passed away in 2019, and I saw Ian at his place two times after that.
The last time was in June of 2023.  When I first arrived, Ian emphatically  told me he hadn't thought about mathematics in 15 years.
He pointed to the Sudoku puzzles and the newspapers on his table as evidence.  We chatted about mutual friends in mathematics
and other memories.  

One evening during my visit, Ian and his sister and I went across the road for dinner -- fish and chips and beer.  During that dinner
it came out in conversation with Ian's sister -- Ian had indeed recently been fiddling with some mathematics, and Ian told me 
about the Clifford circle for the $n$-line.  After dinner, when I was back at his place chatting; at some point, Ian lifted himself out of his chair,
walked over to the other side of the room, picked up a manuscript, and gifted it to me.  He explained that that was what he had been fiddling with and that it was a supplementary chapter to a book he had written just after high school on circles and triangles.  It seems that the manuscript to the book was lost, but I was being given the supplementary chapter.  I didn't quite know what to make of that, but I carefully packed it in my suitcase for my trip home.

After Ian passed away, his son kindly sent me scans of the original handwritten manuscript of the supplementary chapter
and the tex source of the printed copy that Ian gave to me.  Since the topic was lines and circles in the plane, I got a few undergraduates together to work through the manuscript.  

The author of this manuscript was a talented math student right out of high school.  He clearly had not read our key reference for symmetric functions -- his notations for symmetric functions are certainly nonstandard for anyone that has read the Symmetric Functions Bible.  This student shows a  penchant for thorough work and thinking.  For the first main theorem appearing in the manuscript he gives 6 or 7 different proofs, all from different points of view, before moving on to generalizations.  This high school student is incredibly deft with classical and projective geometry and complex numbers (linear equations, determinants, lemniscates, cardiods, deltoids, Euler lines, coaxal systems, Newton identities for symmetric functions, etc).  Some of the induction proofs are a little bit clumsy -- it seems that this student has not been formally taught `proof by induction' like we might do in a first proof course for undergraduates.  The command and thoroughness that this student exhibits extends to his referencing of the literature -- in our modern times most of our community has no idea of the main players of classical intersection geometry any more.  But this high school student was on top of this literature.  If there were one piece of advice that I'd give to this student, it would be to read the books of Ian Macdonald and improve his writing style by emulating the master (admittedly, these books were not yet available).

After getting a feel for the contents of this high school student's manuscript I began to understand Macdonald's early trajectory in mathematics.  He did Tripos at Cambridge and had some exposure to the professors there.  Particularly from the vantage of Hodge, Pedoe and Todd, intersection theory and its connection to cohomology was ``in the air" but not fully developed.  Indeed, in his 1958 paper, Macdonald thanks ``Dr. J.A. Todd for his interest and helpful advice".  By the time of his 1962 paper, Macdonald was clearly following the work of Grothendieck, and had understood that cohomology was an efficient way to compute intersections of the type that he had been computing in high school.
In his 1962 paper he already wields the tools of sheaves and cohomology like a master and writes, ``we obtain natural proofs of the results of an earlier paper [5] which were there obtained laboriously by classical methods."  It is truly amazing to see how this high school student's interest in intersections in classical geometry led him to the very forefront of the technology of cohomology and algebraic geometry that was being vigorously developed at the time and to his proof of an important case of the Weil conjectures.

By 1962, without a Ph.D., Ian Macdonald was no longer a high school student, but had followed his nose to already become a mature mathematician of the highest caliber and a great expositior.

$$
\begin{matrix}
\vcenter{\hbox{\includegraphics[scale=0.65]{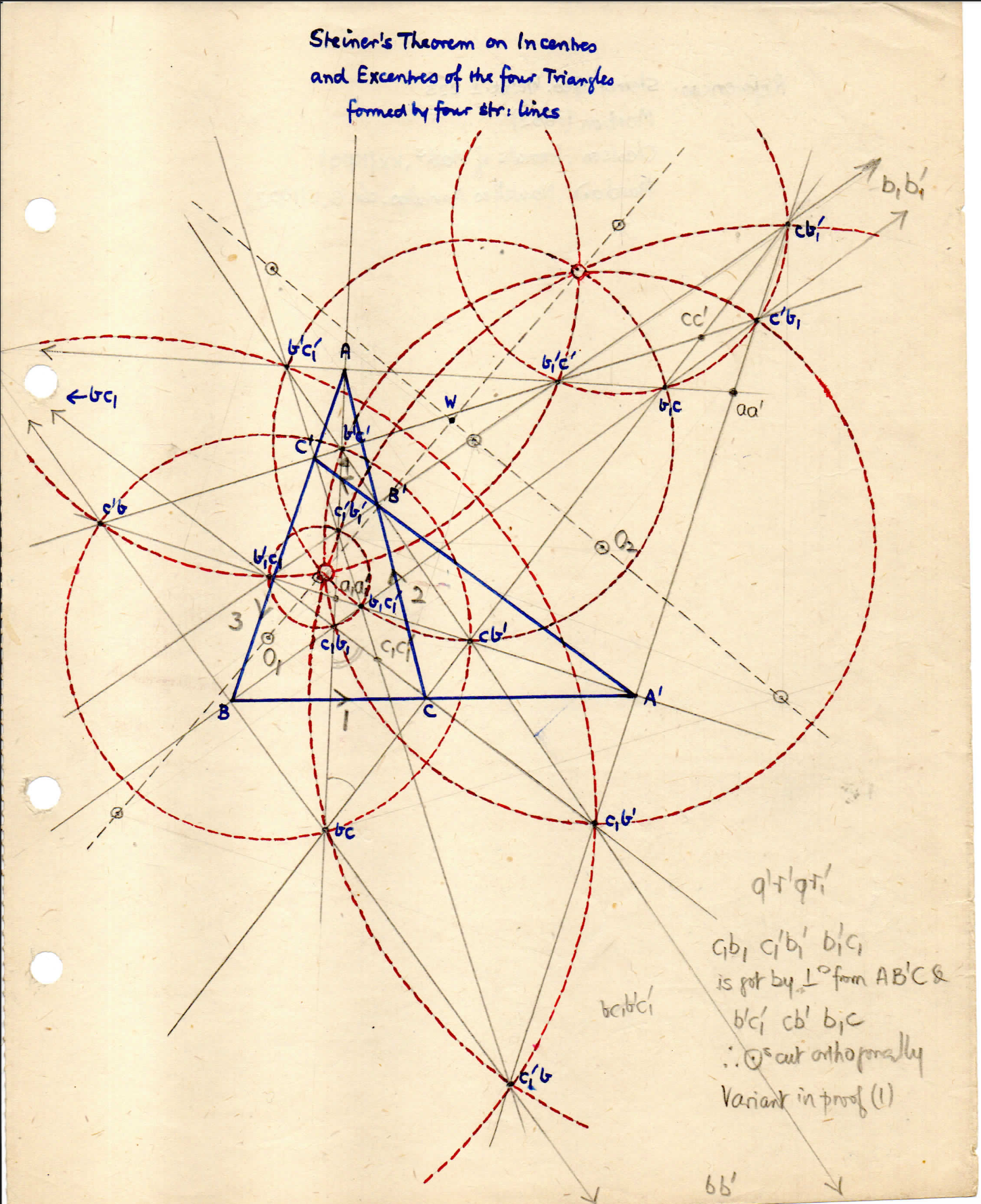}}} \\
\\
\hbox{A diagram from Ian Macdonald's 1947 manuscript on the $n$-line}
\end{matrix}
$$

\subsection{Clifford's $n$-line chain}

\bigskip

Two generic lines $\ell_1$ and $\ell_2$ intersect in a point $A_{12}$.  The point $A_{12}$ is the \emph{Clifford point} of the $2$-line.
$$
\begin{matrix}
\begin{tikzpicture}[scale=0.3]
    \tkzDefPoints{4/2/A,-4/1/B}
    \tkzDrawLine[color = blue](A,B)
    \tkzLabelLine[pos=-0.25](A,B){$\ell_1$}
    \tkzDefPoints{-0.5/6/C,-1.5/0/D}
    \tkzDrawLine[color = blue](C,D)
    \tkzLabelLine[pos=-0.25](C,D){$\ell_2$}
  
    \tkzInterLL(A,B)(C,D)
        \tkzGetPoint{A_{12}}
    \tkzDrawPoint(A_{12})
    \tkzLabelPoint[below left](A_{12}){$A_{12}$}
  \end{tikzpicture}
\end{matrix}
$$
Each pair of lines in a generic $3$-line $\{\ell_1, \ell_2,\ell_3\}$ intersect in a point, and these three points determine a circle $c$.  
The circle $c$ is the \emph{Clifford circle of the $3$-line}.
$$
\begin{matrix}
\begin{tikzpicture}[scale=0.5]
    \tkzDefPoints{4/2/A,-4/1/B}
    \tkzDrawLine[color = blue](A,B)
    \tkzLabelLine[pos=-0.25](A,B){$\ell_1$}
    \tkzDefPoints{-0.5/6/C,-1.5/0/D}
    \tkzDrawLine[color = blue](C,D)
    \tkzLabelLine[pos=-0.25](C,D){$\ell_2$}
    \tkzDefPoints{-1/6/E, 2/0/F}
    \tkzDrawLine[color = blue](E,F)
    \tkzLabelLine[pos=-0.25](E,F){$\ell_3$}
    
    \tkzInterLL(A,B)(C,D)
        \tkzGetPoint{A_{12}}
    \tkzDrawPoint(A_{12})
    \tkzLabelPoint[below left](A_{12}){$A_{12}$}
    \tkzInterLL(A,B)(E,F)
        \tkzGetPoint{A_{13}}
    \tkzDrawPoint(A_{13})
    \tkzLabelPoint[below, yshift=-5, xshift=-3](A_{13}){$A_{13}$}
    \tkzInterLL(C,D)(E,F)
        \tkzGetPoint{A_{23}}
    \tkzDrawPoint(A_{23})
    \tkzLabelPoint[above left, xshift=-2, yshift=-2](A_{23}){$A_{23}$}

    \tkzDefCircle[circum](A_{12},A_{13},A_{23})
        \tkzGetPoint{C_4}
    \tkzDrawCircle[thick, black](C_4,A_{12})    
    \tkzLabelCircle[below](C_4,A_{12})(250){$c$}    
\end{tikzpicture}
\end{matrix}
$$
Each triple of lines in a generic $4$-line $\{\ell_1, \ell_2, \ell_3,, \ell_4\}$
determines a Clifford circle, giving the circles $c_1, c_2, c_3, c_4$.  These four circles intersect in a point $W$.  
The point $W$ is the \emph{Clifford point of the $4$-line}.
$$
\begin{matrix}
\begin{tikzpicture}
    \tkzDefPoints{4/2/A,-4/1/B}
    \tkzDrawLine[color = blue](A,B)
    \tkzLabelLine[pos=-0.25](A,B){$\ell_1$}
    \tkzDefPoints{-0.5/6/C,-1.5/0/D}
    \tkzDrawLine[color = blue](C,D)
    \tkzLabelLine[pos=-0.25](C,D){$\ell_2$}
    \tkzDefPoints{-1/6/E, 2/0/F}
    \tkzDrawLine[color = blue](E,F)
    \tkzLabelLine[pos=-0.25](E,F){$\ell_3$}
    \tkzDefPoints{-4/5/G,4/1.5/H}
    \tkzDrawLine[color = blue](G,H)
    \tkzLabelLine[pos=-0.25](G,H){$\ell_4$}

    \tkzInterLL(A,B)(C,D)
        \tkzGetPoint{A_{12}}
    \tkzDrawPoint(A_{12})
    \tkzLabelPoint[below left](A_{12}){$A_{12}$}
    \tkzInterLL(A,B)(E,F)
        \tkzGetPoint{A_{13}}
    \tkzDrawPoint(A_{13})
    \tkzLabelPoint[below, yshift=-5, xshift=-3](A_{13}){$A_{13}$}
    \tkzInterLL(A,B)(G,H)
        \tkzGetPoint{A_{14}}
    \tkzDrawPoint(A_{14})
    \tkzLabelPoint[below right, xshift=-4, yshift=-3](A_{14}){$A_{14}$}
    \tkzInterLL(C,D)(E,F)
        \tkzGetPoint{A_{23}}
    \tkzDrawPoint(A_{23})
    \tkzLabelPoint[above left, xshift=-2, yshift=-2](A_{23}){$A_{23}$}
    \tkzInterLL(C,D)(G,H)
        \tkzGetPoint{A_{24}}
    \tkzDrawPoint(A_{24})
    \tkzLabelPoint[left, xshift=-1](A_{24}){$A_{24}$}
    \tkzInterLL(E,F)(G,H)
        \tkzGetPoint{A_{34}}
    \tkzDrawPoint(A_{34})
    \tkzLabelPoint[below left, yshift=-7, xshift=2](A_{34}){$A_{34}$}

    \tkzDefCircle[circum](A_{23},A_{24},A_{34})
        \tkzGetPoint{C_1}
    \tkzDrawCircle[thick, black](C_1,A_{23})
    \tkzLabelCircle[above right](C_1,A_{23})(-50){$c_1$}
    \tkzDefCircle[circum](A_{13},A_{14},A_{34})
        \tkzGetPoint{C_2}
    \tkzDrawCircle[thick, black](C_2,A_{13})
    \tkzLabelCircle[right](C_2,A_{13})(180){$c_2$}
    \tkzDefCircle[circum](A_{12},A_{14},A_{24})
        \tkzGetPoint{C_3}
    \tkzDrawCircle[thick, black](C_3,A_{12})
    \tkzLabelCircle[right](C_3,A_{12})(200){$c_3$}
    \tkzDefCircle[circum](A_{12},A_{13},A_{23})
        \tkzGetPoint{C_4}
    \tkzDrawCircle[thick, black](C_4,A_{12})    
    \tkzLabelCircle[above left](C_4,A_{12})(250){$c_4$}

    \tkzInterCC[common=A_{34}](C_1,A_{23})(C_2,A_{13})
        \tkzGetFirstPoint{W}
    \tkzLabelPoint[above right](W){$W$}


    \tkzDrawPoint(W)

\end{tikzpicture}
\end{matrix}
$$
Each $4$-tuple of lines in a generic $5$-line $\{\ell_1, \ell_2, \ell_3,, \ell_4, \ell_5\}$
determines a Clifford point, giving the points $p_1, p_2, p_3, p_4, p_5$.  These five points lie on a circle $C$.
The circle $C$ is the \emph{Clifford circle of the $5$-line}. $\quad\ldots\quad$ and so on $\quad\ldots\quad$

\subsection{Ian Macdonald's general formulation}

Let $y_1,\ldots, y_n\in \CC^\times$. For $i\in \{1, \ldots, n\}$ let $\ell_i$ be the line consisting of the points in $\CC$ that are equidistant 
from $0$ and $y_i$.  The set of $n$ lines is \emph{the $n$-line} $\cL=\{\ell_1, \ldots, \ell_n\}$, where
$$\ell_i =\{ z\in \CC\ |\ \bar z = t_i(z-y_i)\}, \qquad\hbox{where}\quad t_i = \frac{-\overline{y_i}}{y_i}.
\qquad\quad
\begin{matrix}
\begin{tikzpicture}[xscale=2,yscale=1]
\node at (0, 0) {\tiny$\bullet$};  
    \filldraw (0, 0) node[anchor=north,yshift=-0.0cm] {\tiny$0$};
\node at (2, 0) {\tiny$\bullet$};  
    \filldraw (2, 0) node[anchor=north,yshift=-0.0cm] {\tiny$y_i$};
\draw (0,0) -- (2,0);
\draw (1,-1) -- (1,1);
\node at (1.15, .75) {\tiny$\ell_i$};  
\draw (0,0) -- (1,0.3) -- (2,0);
\draw (0,0) -- (1,-0.3) -- (2,0);
\draw (0,0) -- (1,0.5) -- (2,0);
\draw (0,0) -- (1,-0.5) -- (2,0);
%
\end{tikzpicture}
\end{matrix}
$$
For $k\in \{0,1, \ldots, n-1\}$, define
$$c_k(\cL) = \frac{y_1t_1^{n-1-k}}{g_1(\cL)}+\frac{y_2t_2^{n-1-k}}{g_2(\cL)}+\cdots+\frac{y_nt_n^{n-1-k}}{g_n(\cL)},$$
where
$$g_j(\cL) = (t_j-t_1)(t_j-t_2)\cdots (t_j-t_{j-1})(t_j-t_{j+1})(t_j-t_{j+2})\cdots (t_j-t_n),$$
for $j\in \{1, \ldots, n\}$. 


\begin{thm} [Clifford's chain] Let $\cL=\{\ell_1, \ldots, \ell_n\}$ be an $n$-line (satisfying an appropriate genericity condition).
\item[] Case $n$ even: Each $(n-1)$-subset of the $n$-line determines a Clifford circle, and these $n$ Clifford circles intersect in 
a unique point $p(\cL)$.  Let $k\in \ZZ_{>0}$ such that $n = 2k$ and 
let $a_1,\ldots, a_{k-1}\in \CC$ be given by
$$
\begin{pmatrix} a_1 \\ a_2 \\ \vdots \\ a_{k-1} \end{pmatrix}
=
\begin{pmatrix}
c_2(\cL) &\cdots &c_k(\cL) \\
\vdots &&\vdots \\
c_k(\cL) &\cdots &c_{2k-2}(\cL)
\end{pmatrix}^{-1}
\begin{pmatrix} -c_1(\cL) \\ -c_2(\cL)  \\ \vdots \\ -c_{k-1}(\cL) \end{pmatrix}.
$$
Then
$$
p(\cL) 
= c_0(\cL)+a_1 c_1(\cL)+a_2c_2(\cL)+\cdots+a_{k-1}c_{k-1}(\cL)
$$
is the Clifford point of the $n$-line $\cL = \{\ell_1, \ldots, \ell_{2k}\}$.
\item[] Case $n$ odd:  Each $(n-1)$-subset of the $n$-line determines a Clifford point, and these $n$ Clifford points lie on
a unique circle $C(\cL)$.  Let $k\in \ZZ_{>0}$ such that $n = 2k+1$.
The Clifford circle $C(\cL)$ is given by
$$C(\cL) = 
\{ A(\cL) - \theta B(\cL)\ |\ \theta\in U_1(\CC)\},
\qquad\hbox{where}\qquad
U_1(\CC) = \{\theta\in \CC\ |\ \theta\bar \theta = 1\},
$$
$$A(\cL) = \frac{
\det\begin{pmatrix}
c_0(\cL) &\cdots &c_{k-1}(\cL) \\
c_1(\cL) &\cdots &c_k(\cL) \\
\vdots &&\vdots \\
c_{k-1}(\cL) &\cdots &c_{2k-2}(\cL)
\end{pmatrix}
}{
\det\begin{pmatrix}
c_2(\cL) &\cdots &c_k(\cL) \\
\vdots &&\vdots \\
c_k(\cL) &\cdots &c_{2k-2}(\cL)
\end{pmatrix}
}
\qquad\hbox{and}\qquad
B(\cL) = \frac{
\det\begin{pmatrix}
c_1(\cL) &\cdots &c_k(\cL) \\
\vdots &&\vdots \\
c_k(\cL) &\cdots &c_{2k-1}(\cL)
\end{pmatrix}
}{
\det\begin{pmatrix}
c_2(\cL) &\cdots &c_k(\cL) \\
\vdots &&\vdots \\
c_k(\cL) &\cdots &c_{2k-2}(\cL)
\end{pmatrix}
}.$$
\end{thm}


$$
\begin{matrix}
\vcenter{\hbox{\includegraphics[scale=0.5]{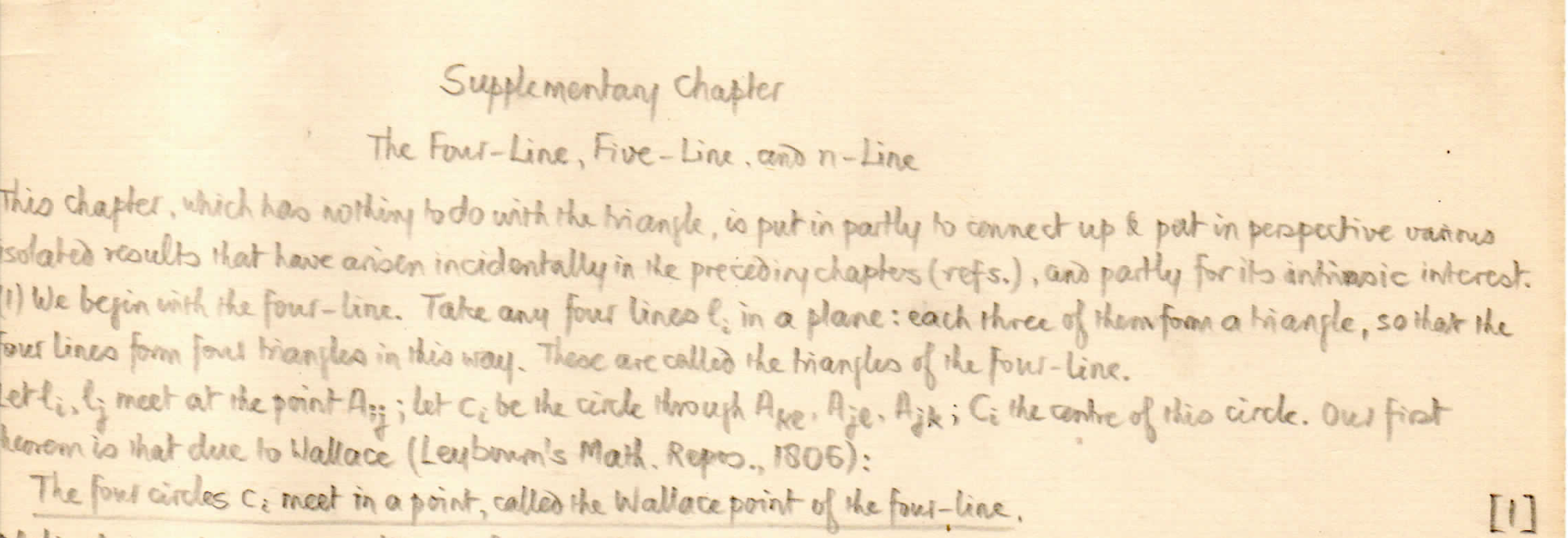}}} \\
\\
\hbox{The first page of the high school student's manuscript}
\end{matrix}
$$

\bigskip\bigskip\bigskip

\section{The symmetric product of a curve $\Sigma$}

After high school Macdonald served in the military and then did the Mathematical Tripos at Trinity.  After finishing at Cambridge in 1952,
at the insistence of his father, Macdonald took competitive exams for a civil service job  (i.e. a government job).
He ``stuck it out for five years" in a civil service job  before leaving his good secure job for 
a temporary (1957--1960) 
position at 
Manchester and then another temporary 
position 
(1960--1963) at Exeter University.  Then he became ``Fellow and Tutor in Mathematics" at Magdalen College at Oxford until 1972.

In 1958, Macdonald's first paper appeared in Proceedings of the Cambridge Philosophical Society.
Very likely this study arose as a continuation of his study of intersections of lines and circles.  The paper is entitled ``Some enumerative formulae for algebraic curves".
In Part I, Macdonald gives a generalization of de Jonqui\`eres formula and Part II makes contact with Schur functions and Schubert conditions in intersection 
theory.  It shows a mastery of the methods of the classical Italian algebraic geometry school.  This paper is a significant
development of his high school knowledge of intersection theory.  Even so, it hardly gives any hint of the amazing achievement that
was to come next.

By 1962, Macdonald had understood that intersection numbers of families of ourves could be computed by using cohomology
as a tool. 
In his paper on the cohomology of symmetric products of an algebraic curve \cite{Mac62b} he states ``in particular we obtain natural proofs of the results of an earlier paper \cite{Mac58} which were there
obtained laboriously by classical methods.''

\newpage

\subsection{The cohomology of a symmetric product of a curve}

Let $\Sigma$ be a curve.  
$$
\begin{matrix}
\begin{tikzpicture}
\draw[smooth] (0,1) to[out=30,in=150] (2,1) to[out=-30,in=210] (3,1) to[out=30,in=150] (5,1) to[out=-30,in=210] (6,1) to[out=30,in=150] (8,1) to[out=-30,in=30] (8,-1) to[out=210,in=-30] (6,-1) to[out=150,in=30] (5,-1) to[out=210,in=-30] (3,-1) to[out=150,in=30] (2,-1) to[out=210,in=-30] (0,-1) to[out=150,in=-150] (0,1);
\draw[smooth] (0.4,0.1) .. controls (0.8,-0.25) and (1.2,-0.25) .. (1.6,0.1);
\draw[smooth] (0.5,0) .. controls (0.8,0.2) and (1.2,0.2) .. (1.5,0);
\draw[smooth] (3.4,0.1) .. controls (3.8,-0.25) and (4.2,-0.25) .. (4.6,0.1);
\draw[smooth] (3.5,0) .. controls (3.8,0.2) and (4.2,0.2) .. (4.5,0);
\draw[smooth] (6.4,0.1) .. controls (6.8,-0.25) and (7.2,-0.25) .. (7.6,0.1);
\draw[smooth] (6.5,0) .. controls (6.8,0.2) and (7.2,0.2) .. (7.5,0);
\end{tikzpicture}
\end{matrix}
$$
The $n$th symmetric power of $\Sigma$ is 
$$\Sigma(n) = \Sigma^n/S_n,
\qquad\hbox{where}\quad
w\cdot (p_1,\ldots, p_n) = (p_{w^{-1}(1)}, \ldots, p_{w^{-1}(n)}),
$$
for $w\in S_n$ and $(p_1, \ldots, p_n)\in \Sigma^n$.

Cohomology is a creature (more precisely, a functor) that eats spaces and outputs graded rings.
$$
\hbox{spaces} \qquad\qquad
\begin{matrix}
\begin{tikzpicture}[scale=2]
\coordinate (A) at (0,0);
\coordinate (B) at (-150:1cm);
\coordinate (C) at (150:1cm);
\draw (B) -- (A) -- (C) ;
\draw (-150:1cm) arc [start angle=-150, end angle=-10, radius=1cm];
\draw (150:1cm) arc [start angle= 150, end angle=10, radius=1cm];
\draw (B) -- (-160:0.9cm) -- (-150:0.8cm) -- (-160:0.7cm) -- (-150:0.6cm) -- (-160:0.5cm) -- (-150:0.4cm) -- (-160:0.3cm) -- (-150:0.2cm);
\draw (C) -- (160:0.9cm) -- (150:0.8cm) -- (160:0.7cm) -- (150:0.6cm) -- (160:0.5cm) -- (150:0.4cm) -- (160:0.3cm) -- (150:0.2cm);
\draw (10:1cm) arc [start angle = 90, end angle = 270, x radius=2pt, y radius= 5pt];
\draw (10:1cm) -- (7:1.2cm);
\draw (-10:1cm) -- (-7:1.2cm);
\draw (1.2,0) ellipse [x radius=1pt, y radius= 4pt];
\draw (0,0.5) node{$H^*$};
\draw (-0.7,0) node{$\Sigma(n)$};
\draw (1.8,0) node{$H^*(\Sigma(n),\ZZ)$};
\end{tikzpicture}
\end{matrix}
\qquad\qquad\hbox{graded rings}
$$
In spite of the frightening teeth, cohomologies are really quite friendly (it is the spaces that are dangerous).
How does one write down the cohomology $H^*(X;\ZZ)$ of a space $X$?  Well, $H^*(X,\ZZ)$ is a graded ring
and a graded ring is written down in a presentation by generators and relations.  Macdonald's 1962 paper
gives an elegant presentation of the graded ring $H^*(\Sigma(n),\ZZ)$, the cohomology of the $n$th symmetric product of
a curve $\Sigma$.

\begin{thm} Let $\Sigma$ be a curve of genus $g$.  The cohomology ring $H^*(\Sigma(n),\ZZ)$ is the
$\ZZ$-algebra presented by generators
$$\xi_1, \ldots, \xi_g, \quad \xi_1', \ldots, \xi_g', \quad \eta$$
and relations
\item[(a)] If $i,j\in \{1, \ldots, g\}$ then
$$\xi_i\xi_j = -\xi_j\xi_i, \qquad \xi_i'\xi_j' = -\xi_j' \xi_i', \qquad \xi_i \xi_j' = -\xi_j' \xi_i, $$
$$\xi_i\eta = \eta\xi_i, \qquad\qquad \xi_i' \eta = \eta\xi_i',$$
\item[(b)] If $a,b,c,q\in \ZZ_{\ge 0}$ and $a+b+2c+q=n+1$ and $i_1, \ldots, i_a, j_1, \ldots, j_b, k_1, \ldots, k_c$ are distinct elements of 
$\{1, \ldots, g\}$ then
$$\xi_{i_1}\cdots \xi_{i_a} \xi'_{j_1}\cdots \xi'_{j_b}(\xi_{k_1}\xi'_{k_1}-\eta)\cdots (\xi_{k_c}\xi'_{k_c}-\eta)\eta^q = 0.$$
\end{thm}



\subsection{The Weil conjectures for the symmetric product $\Sigma(n)$ of a curve}

Weil's famous conjectures about zeta functions of algebraic varieties are from his paper of 1949 \cite{We49}. These conjectures
stimulated a huge effort which included the development of \'etale cohomology and $\ell$-adic cohomology.  The Weil conjectures were proved in the 1960s and 70s: the proof of the rationality conjecture came in 1960 (Dwork), 
the proof of the functional equation and Betti numbers connection in 1965 (Grothendieck school) and the
analogue of the Riemann hypothesis in 1974 (Deligne).  
In 1962, Macdonald proved Weil's conjectures in an important special case: ``... we calculate the zeta function of 
$\Sigma(n)$ and verify Weil's conjectures in this case."

The zeta function $Z(t)$ of an algebraic variety $X$ is an exponential generating function for the number of points of $X$ over the
finite fields $\FF_{q^n}$,
$$\frac{d\ }{dt} \log Z(t) = \sum_{n\in \ZZ_{>0}} \Card(X(\FF_{q^n})) t^{n-1}.$$
Let $\Sigma$ be a curve of genus $g$ and assume that $\rho_1, \ldots, \rho_{2g}\in \CC$ are such that 
$$Z_1(t) = \frac{(1-\rho_1t)\cdots (1-\rho_{2g}t)}{(1-t)(1-qt)}
\qquad
\hbox{is the zeta function of $\Sigma$.}
$$
Let $\phi_0(t) = 1-t$ and, for $k\in \{1,\ldots, 2g\}$, let
$$\phi_k(t) = \prod_{1\le i_1<\cdots < i_k \le 2g} (1-\rho_{i_1}\cdots \rho_{i_k}t).$$
Then define
$$
F_k(t) = \begin{cases}
\phi_k(t)\phi_{k-2}(t)\phi_{k-4}(t) \cdots, &\hbox{if $k\in \{0, 1, \ldots, n\}$,} \\
F_{2n-k}(q^{k-n}t), &\hbox{if $k\in \{n+1, \ldots, 2n\}$.}
\end{cases}
$$

\begin{cor} The Weil conjectures hold for $\Sigma(n)$.  More specifically,
\item[(a)] The zeta function of $\Sigma(n)$ is 
$$Z_n(t) = \frac{F_1(t)F_3(t)\cdots F_{2n-1}(t)}{F_0(t)F_2(t)\cdots F_{2n}(t)}.$$
\item[(b)] The Riemann hypothesis for $\Sigma(n)$ holds: 
$$\hbox{All roots of $Z_n(t)$ have absolute value in }\quad
\{ q^{-\frac12\cdot 0}, q^{-\frac12\cdot 1}, q^{-\frac12\cdot 2}, \ldots, q^{-\frac12\cdot 2n}\}.$$
\item[(c)] The functional equation for $\Sigma(n)$ is
$$Z_n\big(\hbox{$\frac{1}{q^nt}$}\big) = (-q^{-\frac12n}t)^{(-1)^n\binom{2g-2}{n}} Z_n(t).$$
\end{cor}

%

\newpage

\section{I.G.\ Macdonald as influencer}

\subsection{Deligne-Lusztig 1976}

In Lecture Notes in Math.\ \textbf{131}, T.\ Springer precisely states 
conjectures of Macdonald about complex representations of finite groups of Lie type.
Looking back at these references, one gathers that the notes of Macdonald on Hall polynomials 
that were circulating in the late 1960's eventuallly became Chapter IV of his book on Symmetric functions and Hall polynomials.  
T.\ Springer's expositions appearing in \cite{Spr70} make it clear  that, by 1968,  Ian Macdonald had understood how the 
type $GL_n$ story from J.A. Green's 1955 paper could be reshaped for a statement for general Lie types.  
Macdonald's conjectures were proved by Deligne and Lusztig in 1976.


\bigskip\bigskip

$$
\begin{matrix}
\vcenter{\hbox{\includegraphics[scale=0.6]{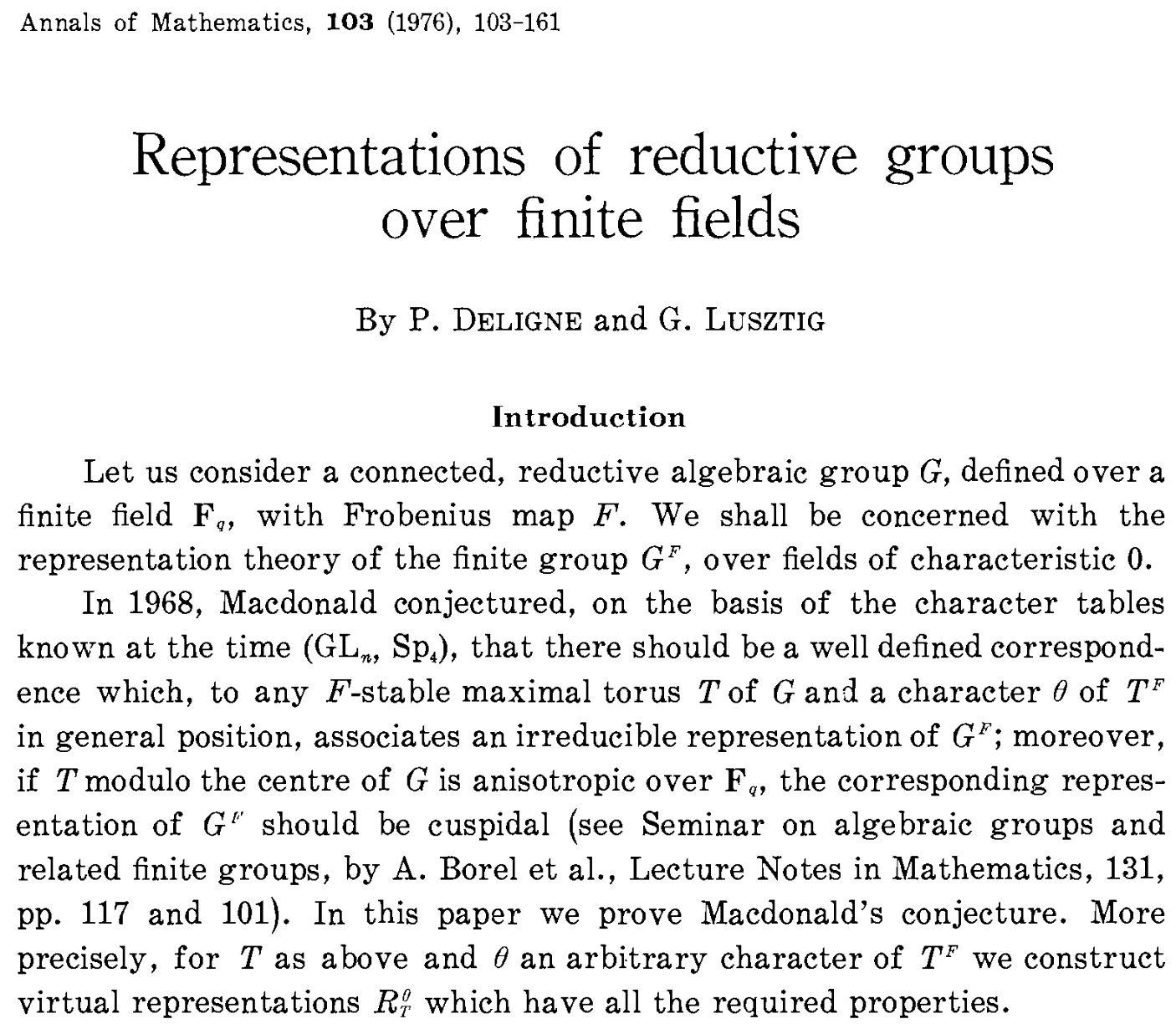}}} \\
\end{matrix}
$$

\newpage

\subsection{Maulik-Yun 2013}

It is still the case that most topologists and geometers view Macdonald's computation of the cohomology of the symmetric product of a curve
\cite{Mac62b} as his most well known achievement.  
In recent years the study of moduli spaces of curves and related cohomological Hall algebras has become an important part of
geometry and mathematical physics, and Macdonald's study of symmetric products of curves 
continues to be an important stimulus for research in this direction today.

$$
\begin{matrix}
\vcenter{\hbox{\includegraphics[scale=0.5]{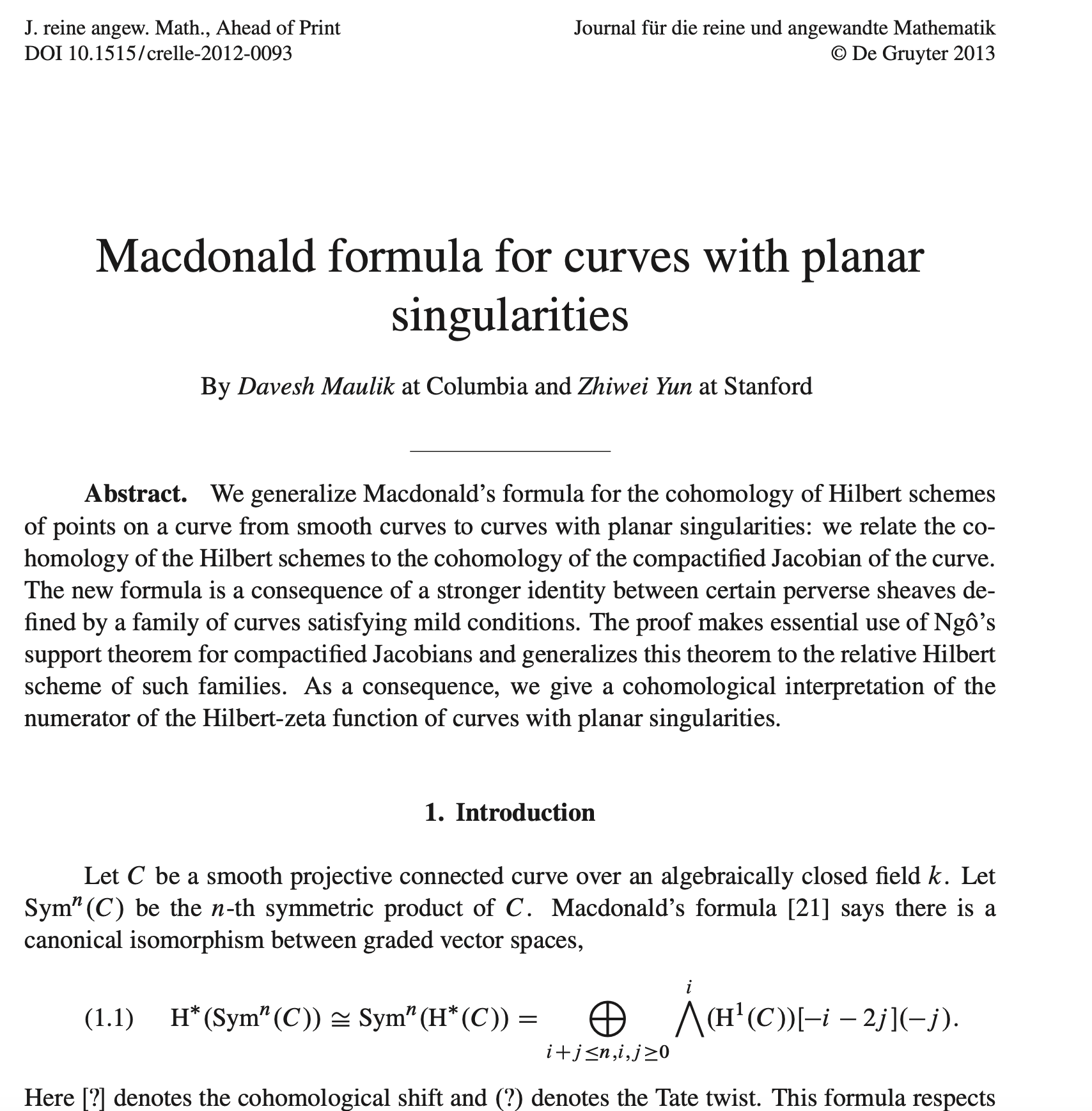}}} \\
\end{matrix}
$$


\newpage

\subsection{Casselman 2012}

An announcement of Macdonald's computation of the spherical function for $p$-adic groups appeared in 1968, and the full
details appeared in his book published by the University of Madras in 1971.  From the point of view of symmetric function theory,
Macdonald proved that the favorite formula \cite[Ch.\ III (2.1)]{Mac} for the Hall-Littlewood polynomial
$$P_\lambda(x;t) = \frac{1}{v_\lambda(t)} \sum_{w\in S_n} w\Big(x_1^{\lambda_1}\cdots x_n^{\lambda_n}\prod_{i<j} \frac{x_i-tx_j}{x_i-x_j}\Big)$$
generalizes to all Lie types and is a formula for the spherical function for $G/K$ where $G$ is the corresponding $p$-adic group $G=G(\QQ_p)$ and $K=G(\ZZ_p)$ is a maximal compact subgroup of $G$.

$$
\begin{matrix}
\vcenter{\hbox{\includegraphics[scale=0.55]{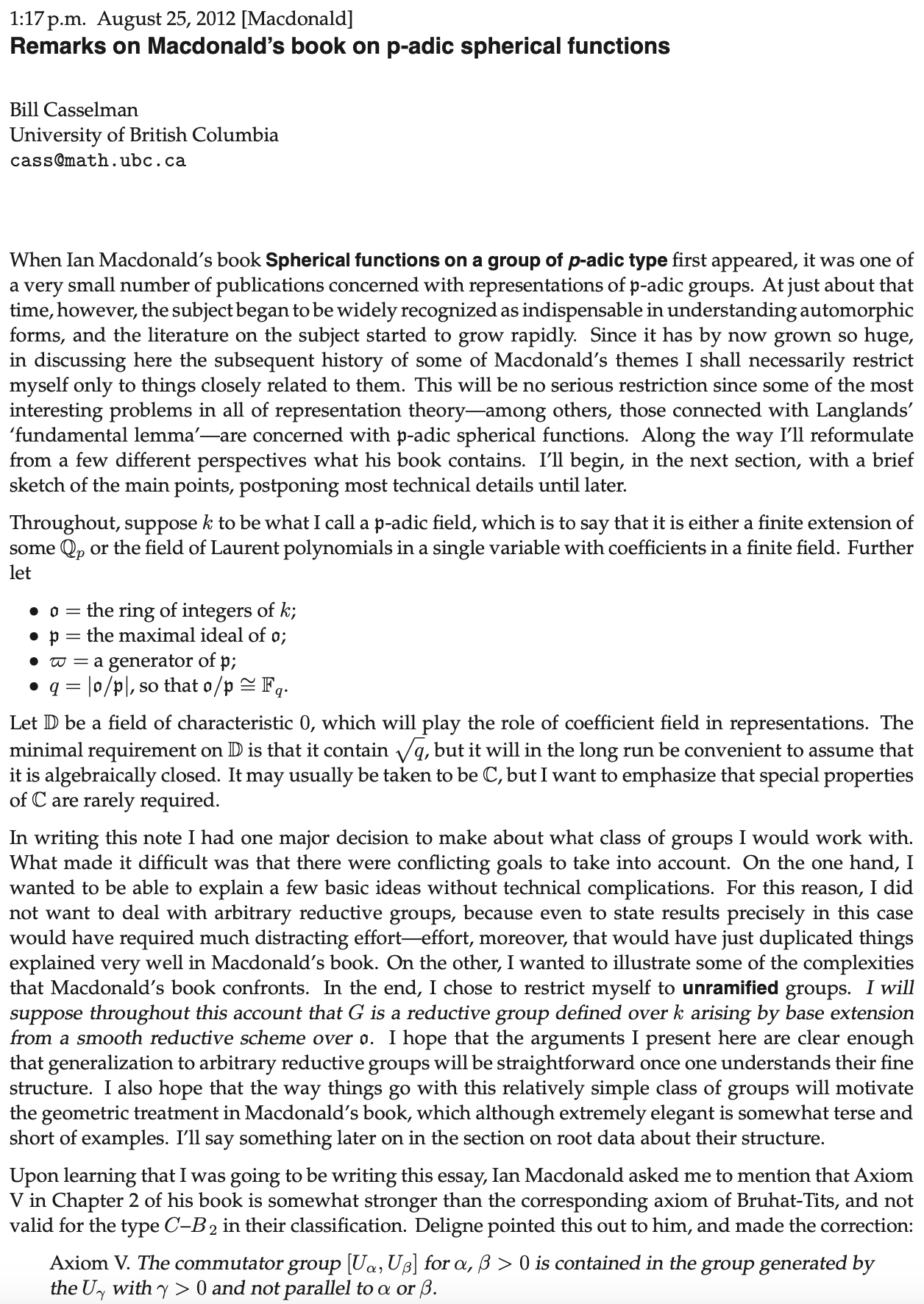}}} 
\end{matrix}
$$

\newpage

\subsection{V. Kac, \emph{Infinite dimensional Lie algebras}, Cambridge University Press, 1982}

Macdonald's work on $p$-adic groups drew him into the combinatorics of affine root systems and he made a thorough classification
and study of affine root systems and affine Weyl groups, resulting in his 1972 paper entitled 
``Affine root systems and Dedekind's $\eta$-function''.
This study brought him into contact with affine Kac-Moody Lie algebras and formulas for characters of their representations.

\bigskip

$$
\begin{matrix}
\vcenter{\hbox{\includegraphics[scale=0.6]{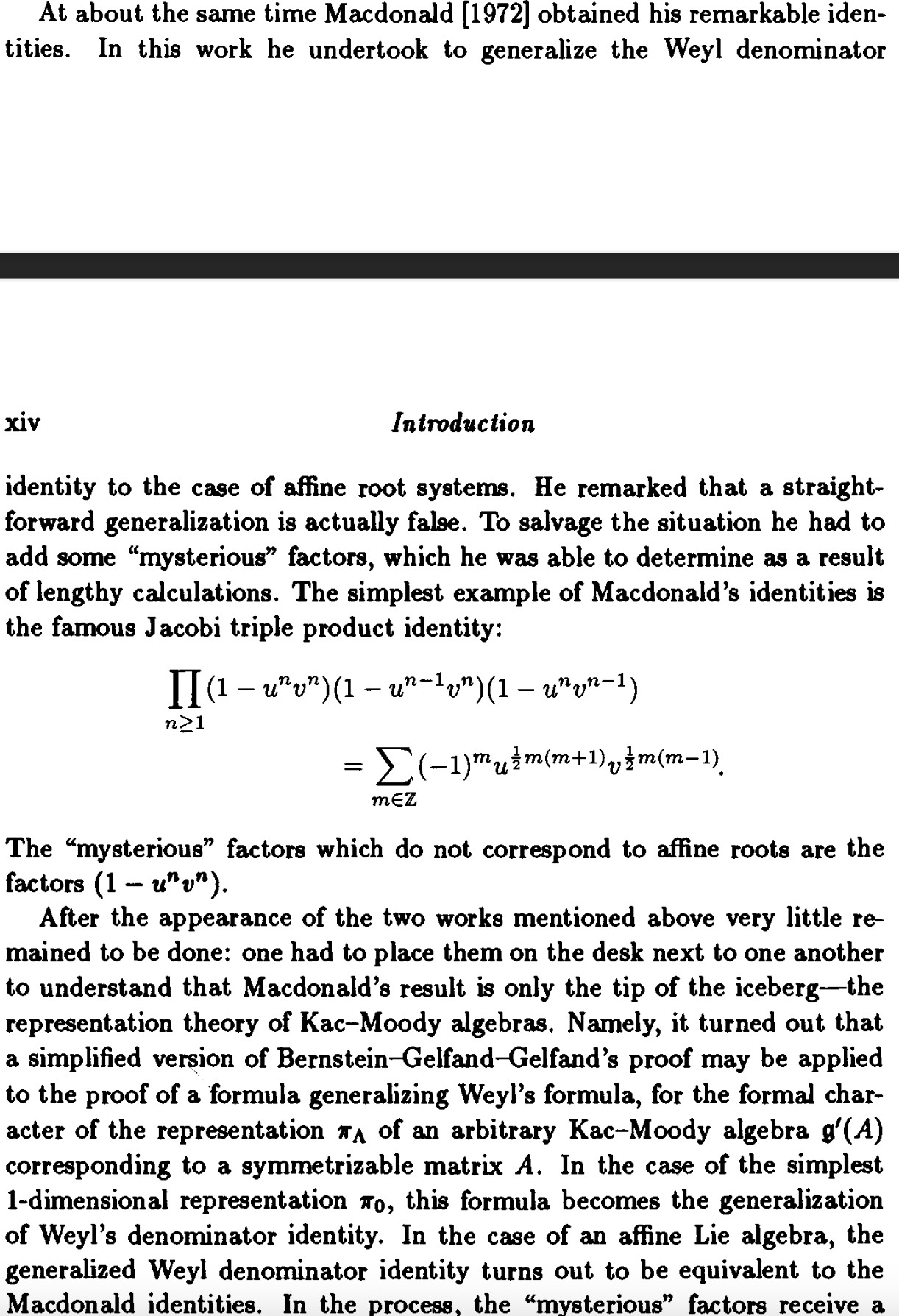}}} 
\end{matrix}
$$

\newpage

\section{I.G.\ Macdonald as translator}

One of Ian Macdonald's great silent contributions to the mathematical community was his work as a translator.

\subsection{I.G.\ Macdonald as translator: Bourbaki}

I.G.\ Macdonald was the first translator of Bourbaki into English.  
It is not clear how much of Bourbaki Macdonald translated as the publisher did not list the 
translator in the published English versions.   A best guess is that the volumes which appeared in English between 1966 and 1974 were translated
by Macdonald.  These volumes comprise more than 2500 pages.

\bigskip

{\large
\begin{enumerate}
\item[] \textbf{Bourbaki, General Topology Parts I and II 1966, vii+437 pp. and iv+363 pp.}
\item[] \textbf{Bourbaki, Theory of Sets 1968, viii+414 pp.}
\item[] \textbf{Bourbaki, Commutative Algebra 1972, xxiv+625 pp.}
\item[] \textbf{Bourbaki, Algebra 1974, xxiii+709 pp.}
\end{enumerate}

\bigskip\noindent


\subsection{I.G.\ Macdonald as translator: Dieudonn\'e}

I.G. Macdonald's work as a translator of Dieudonn\'e's Treatise on Analysis is documented in \cite{Mar} and he is
explicitly listed as translator in the English version of Dieudonn\'e's Panorama of Pure Mathematics, which appeared in 1982.
Together, these volumes amount to more than 2300 pages.

\bigskip

{\large
\begin{enumerate}
\item[] \textbf{Dieudonn\'e, Foundations of Modern Analysis 1960 and 1969, xiv+361 pp.}
\item[] \textbf{Dieudonn\'e, Treatise on Analysis Vol. II 1970 and 1976, xviii+387 pp.}
\item[] \textbf{Dieudonn\'e, Treatise on Analysis Vol. III 1972, xvii+388 pp.}
\item[] \textbf{Dieudonn\'e, Treatise on Analysis Vol. IV 1974, xiv+444 pp.}
\item[] \textbf{Dieudonn\'e, Treatise on Analysis Vol. V 1977, xiv+243 pp.}
\item[] \textbf{Dieudonn\'e, Treatise on Analysis Vol. VI 1978, xi+239 pp}
\item[] \textbf{Dieudonn\'e, A panorama of pure mathematics 1982, x+289 pp.}
\end{enumerate}
}

%




\newpage

\section{I.G. Macdonald for my students}

Every once in a while, not infrequently, a student comes by my office and says ``I'd like to learn about \textbf{Lie groups},
do you have a reference that you can recommend?''  I usually find myself saying, ``How about the notes of Macdonald?''
\begin{enumerate}
\item[] Algebraic structure of Lie groups, Cambridge University Press, 1980. \hfil\break
\href{https://doi.org/10.1017/CBO9780511662683.005}{https://doi.org/10.1017/CBO9780511662683.005 }
\end{enumerate}
Every once in a while, not infrequently, a student comes by my office and says ``I'd like to learn about \textbf{algebraic groups},
do you have a reference that you can recommend?''  I usually find myself saying, ``How about the notes of Macdonald?''
\begin{enumerate}
\item[] Linear algebraic groups, \hfil\break
in Lectures on Lie Groups and Lie Algebras, Cambridge University Press 1995
 \hfil\break
\href{https://doi.org/10.1017/CBO9781139172882}{https://doi.org/10.1017/CBO9781139172882}
\end{enumerate}
Every once in a while, not infrequently, a student comes by my office and says ``I'd like to learn about \textbf{reflection groups},
do you have a reference that you can recommend?''  I usually find myself saying, ``How about the notes of Macdonald?''
\begin{enumerate}
\item[] Reflection groups, unpublished notes 1991. \hfil\break
Available at \href{http://math.soimeme.org/~arunram/resources.html}{http://math.soimeme.org/{\raise.17ex\hbox{$\scriptstyle\sim$}}arunram/resources.html}
\end{enumerate}
Every once in a while, not infrequently, a student comes by my office and says ``I'd like to learn about \textbf{algebraic geometry},
do you have a reference that you can recommend?''  I usually find myself saying, ``How about the book of Macdonald?''
\begin{enumerate}
\item[] Algebraic Geometry - Introduction to schemes, published by W.A. Benjamin 1968. \hfil\break
Available at \href{http://math.soimeme.org/~arunram/resources.html}{http://math.soimeme.org/{\raise.17ex\hbox{$\scriptstyle\sim$}}arunram/resources.html}
\end{enumerate}
Every once in a while, not infrequently, a student comes by my office and says ``I'd like to learn about \textbf{Haar measure, spherical functions and harmonic analysis},
do you have a reference that you can recommend?''  I usually find myself saying, ``How about the book of Macdonald?''
\begin{enumerate}
\item[] Spherical functions on a group of $p$-adic type, University of Madras 1971. \hfil\break
Available at \href{http://math.soimeme.org/~arunram/resources.html}{http://math.soimeme.org/{\raise.17ex\hbox{$\scriptstyle\sim$}}arunram/resources.html}
\end{enumerate}
Every once in a while, not infrequently, a student comes by my office and says ``I'd like to learn about \textbf{Kac-Moody Lie algebras},
do you have a reference that you can recommend?''  I usually find myself saying, ``How about the notes of Macdonald?''
\begin{enumerate}
\item[] Kac-Moody Lie algebras, unpublished notes 1983. \hfil\break
Available at \href{http://math.soimeme.org/~arunram/resources.html}{http://math.soimeme.org/{\raise.17ex\hbox{$\scriptstyle\sim$}}arunram/resources.html}
\end{enumerate}
Every once in a while, not infrequently, a student comes by my office and says ``I'd like to learn about \textbf{flag varieties and
Schubert varieties},
do you have a reference that you can recommend?''  I usually find myself saying, ``How about the notes of Macdonald?''
\begin{enumerate}
\item[] Notes on Schubert polynomials: Appendix: Schubert varieties. \hfil\break
Published by LACIM 1991. \hfil\break
Available at \href{http://math.soimeme.org/~arunram/resources.html}{http://math.soimeme.org/{\raise.17ex\hbox{$\scriptstyle\sim$}}arunram/resources.html}
\end{enumerate}

\newpage

\section{I.G.\  Macdonald as an author of books}

{\sl ``If you see a gap in the literature, write a book to fill it." -- I.G. Macdonald}

{\large
\begin{enumerate}
\item[] \textbf{Atiyah-Macdonald, Introduction to commutative algebra 1969}
\item[] \textbf{Spherical functions on a group of p-adic type 1971}
\item[] \textbf{Symmetric functions and Hall polynomials First Edition 1979}
\item[] \textbf{Kac-Moody Lie algebras: unpublished notes 1983}
\item[] \textbf{Hypergeometric functions: unpublished notes 1987}
\item[] \textbf{Reflection groups: unpublished notes 1991}
\item[] \textbf{Schubert polynomials 1991}
\item[] \textbf{Symmetric functions and Hall polynomials Second Edition 1995}
\item[] \textbf{Linear algebraic groups: in Lectures on Lie groups and Lie algebras 1995}
\item[] \textbf{Affine Hecke algebras and orthogonal polynomials 2003}
\end{enumerate}
}


\bigskip\bigskip

\subsection{The first book: Algebraic Geometry - Introduction to Schemes 1968}


$$
\begin{matrix}
\vcenter{\hbox{\includegraphics[scale=0.5]{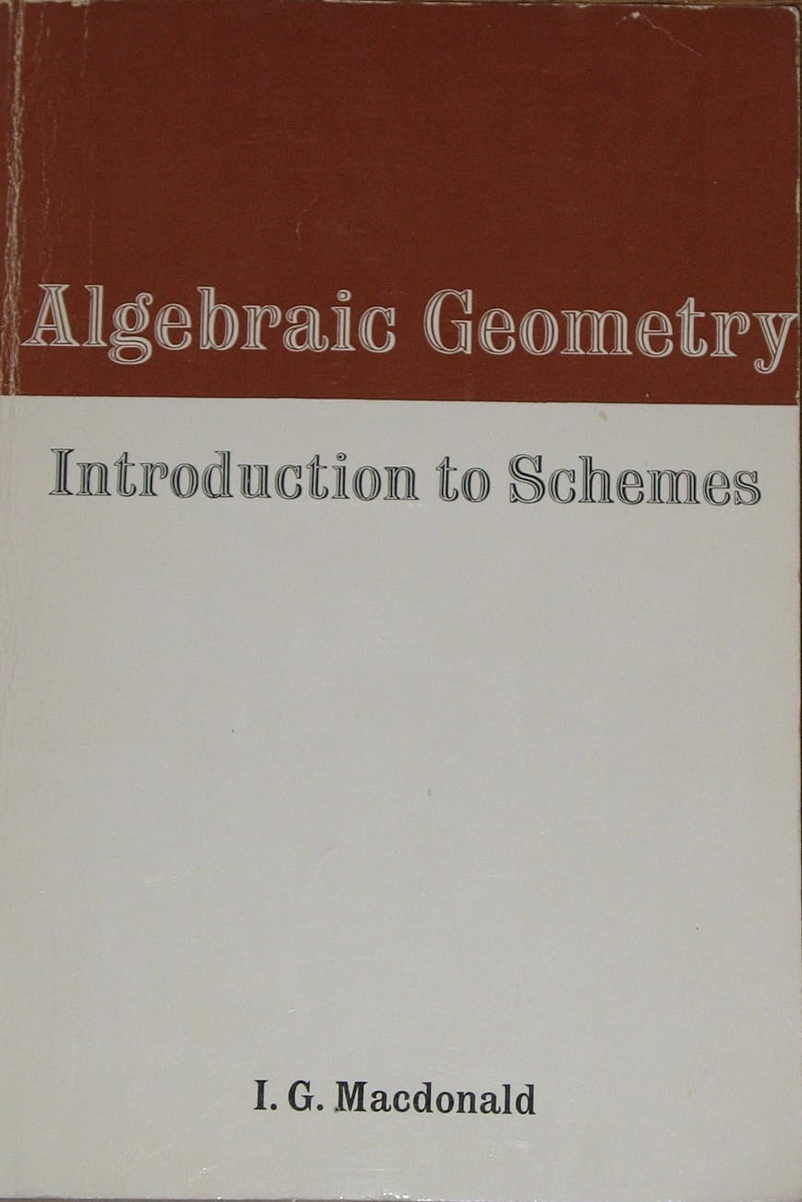}}} \\
\end{matrix}
$$

\newpage

\setcounter{section}{10}
\setcounter{subsection}{0}

\end{document}